\documentclass[11pt,reqno]{amsart}
\usepackage[letterpaper,margin=1in]{geometry}
\usepackage[T1]{fontenc}
\usepackage{lmodern,amsmath,amssymb,amsthm,mathrsfs,microtype}
\usepackage[hidelinks]{hyperref}
\usepackage{enumitem}
\numberwithin{equation}{section}
\newtheorem{theorem}{Theorem}[section]
\newtheorem{proposition}[theorem]{Proposition}
\newtheorem{lemma}[theorem]{Lemma}
\newtheorem{corollary}[theorem]{Corollary}
\theoremstyle{definition}

\theoremstyle{remark}
\newtheorem{remark}[theorem]{Remark}
\newcommand{\Q}{\mathbb Q}
\newcommand{\BC}{\operatorname{BC}}
\title[Functoriality probes over number fields]{Functoriality probes over number fields by mixed zero correlations}
\author{Yangbo Ye}
\address{Department of Mathematics, University of Iowa, Iowa City, IA 52242, USA}
\email{yangbo-ye@uiowa.edu}
\date{Expanded version, October 24, 2026}
\begin{document}
\begin{abstract}
We use an automorphic $L$-function over one number field as a probe for the functorial structure of another automorphic $L$-function over a possibly different field. No automorphic induction is assumed. We first determine the mixed Euler factors arising from the explicit formula, including the residue-degree multiplicities that occur when local traces are grouped by absolute norm. This local analysis gives four probe constructions. Rational-field probes with known base change detect constituent multiplicities and, in equal rank, exact Satake equalities. Probes over a common Galois field detect conjugacy and stabilizers. Dirichlet-character probes give congruence-weighted combinations of Satake traces even when the leading overlap vanishes. A tensor-product construction realizes general different-field probes whenever the required base changes exist. We then separate the pair calculation from the higher connected calculations and work out the cubic and quartic arithmetic explicitly. In an explicit small Fourier box these expressions yield all finite secondary orders of the centered statistics; the required pair prime number theorems also give raw and distinct-index correlations. The main analytic tools are an independent fixed-data proof of Hypothesis~H over number fields, a positive-margin estimate for rationally grouped coefficients, and the theorem of Harcos--Thorner for the realized Rankin--Selberg factors. This expanded version includes the full computations, local examples, and an appendix explaining the character-extension issue avoided by the proof of Hypothesis~H. The resulting probes give concrete tests of functorial relationships, while recovery of automorphic induction remains a longer-term goal.
\end{abstract}
\maketitle
\tableofcontents

\section{Introduction}
Let $E$ and $F$ be number fields, and let $\Pi$ and $\pi$ be unitary cuspidal automorphic representations of $\mathrm{GL}_m(\mathbb A_E)$ and $\mathrm{GL}_n(\mathbb A_F)$, respectively. We regard $L_F(s,\pi)$ as a probe and ask what the mixed zero correlations of
\begin{equation}L_E(s,\Pi)L_F(s,\pi)\label{in:product}\end{equation}
reveal about the functorial structure of $\Pi$. The probe may be defined over $\mathbb Q$, over $E$, or over another number field. Our approach is to identify the arithmetic expressions measured by these correlations and then determine when they detect a base-change constituent, a Galois conjugacy, or a relation among local Satake parameters. We do not assume automorphic induction of $\Pi$ or a factorization of $L_E(s,\Pi)$ over $\mathbb Q$.

Automorphic induction is instead a long-term motivation. At an unramified rational prime $p$, the local Satake matrices of $\Pi$ give the polynomial
\begin{equation}P_{\Pi,p}(X)=\prod_{v\mid p}\det(1-A_{\Pi,v}X^{f_v}),\qquad Nv=p^{f_v}.\label{in:polynomial}\end{equation}
If an automorphic induction exists, these are its standard Euler polynomials. Even without that global transfer, they are well-defined local objects. A suitable family of probes may furnish tests for their compatibility with a global automorphic representation. The results here establish several such tests; they do not assert that the resulting correlation data determine an automorphic induction.

\subsection{Local arithmetic and probe constructions}
The local computation is essential because the explicit formula groups places by their absolute norms. Put $d_E=[E:\mathbb Q]$ and define
\begin{equation}b_\Pi(p^k)=\sum_{\substack{v\mid p\\f_v\mid k}}f_v\operatorname{tr}(A_{\Pi,v}^{k/f_v}).\label{in:b}\end{equation}
These are the power sums of the reciprocal roots of \eqref{in:polynomial}. For the mixed logarithmic series
\begin{equation}K_{\Pi,\pi}(s)=\sum_{p,k\ge1}\frac{(\log p)b_\Pi(p^k)\overline{b_\pi(p^k)}}{p^{ks}},\label{in:K}\end{equation}
Proposition~\ref{id:local-factor-prop} gives exact Euler factors. A pair of residue degrees $f_v,g_w$ contributes exponents involving $\operatorname{lcm}(f_v,g_w)$ with multiplicity $\gcd(f_v,g_w)$. More precisely, writing $L_{v,w}=\operatorname{lcm}(f_v,g_w)$ and $h_{v,w}=\gcd(f_v,g_w)$, the local factor is
\begin{equation}\label{in:local-factor}
\mathscr L_{\Pi,\pi,p}(s)=
\prod_{v\mid p}\prod_{w\mid p}\prod_{i=1}^{m}\prod_{j=1}^{n}
\left(1-\alpha_{v,i}^{L_{v,w}/f_v}
\overline{\beta_{w,j}^{L_{v,w}/g_w}}p^{-L_{v,w}s}\right)^{-h_{v,w}}.
\end{equation}
This accounts for all pairs of places of equal absolute norm, without a Galois assumption. It also explains why merely pairing places along a specified inclusion $F\subset E$ generally misses part of the mixed correlation. The tensor-product identity below supplies the precise global counterpart of these local multiplicities.

We obtain the following four probe constructions, the last of which includes a general realization of the first two.
\begin{enumerate}[label=\textup{(\arabic*)},leftmargin=*]
\item \emph{Rational-field probes with known base change.} Suppose $F=\mathbb Q$ and the specified base change $\tau=\BC_{E/\mathbb Q}\pi$ exists as a unitary isobaric representation. Away from finitely many primes,
\begin{equation}K_{\Pi,\pi}^S(s)=-\frac d{ds}\log L_E^S(s,\Pi\times\widetilde\tau).\label{in:probe}\end{equation}
The residue at $1$ counts the copies of $\Pi$ in $\tau$. If $m=n$, a positive residue is equivalent to $\tau\simeq\Pi$, and hence to the compatible local equalities
\begin{equation}\{\alpha_{\Pi,v,i}\}_{i=1}^m=\{\beta_{\pi,p,i}^{f_v}\}_{i=1}^m\qquad(v\mid p,\ p\notin S),\label{in:Satake}\end{equation}
with the converse supplied by strong multiplicity one. For unequal ranks the residue still detects a global constituent, a conclusion stronger than pointwise divisibility of Euler polynomials.
\item \emph{Probes over a common Galois field.} For $E/\mathbb Q$ Galois with group $G$, and $\Pi,\Sigma$ over $E$, we prove
\begin{equation}K_{\Pi,\Sigma}^S(s)=-\frac d{ds}\log\prod_{g\in G}L_E^S(s,\Pi\times\widetilde{\Sigma^g}).\label{in:Galois}\end{equation}
The residue counts the $g$ for which $\Pi\simeq\Sigma^g$. In particular, the leading self pair term detects $\#\operatorname{Stab}_G(\Pi)$, without assuming induction or nontrivial base change.
\item \emph{Character probes.} A Dirichlet character $\chi$ always has the norm pullback $\chi_E=\chi\circ N_{E/\mathbb Q}$. Its mixed logarithmic series is the negative logarithmic derivative of $L_E(s,\Pi\otimes\overline{\chi_E})$. For $m>1$ the leading mixed overlap is zero, but the secondary pair constants remain arithmetic invariants. Varying $\chi$ and applying character orthogonality gives congruence-weighted combinations of the grouped Satake traces, as in \eqref{pr:orthogonality}.
\item \emph{General different-field probes.} Write $E\otimes_{\mathbb Q}F=\prod_\nu K_\nu$. If the specified unitary isobaric base changes of both representations to every $K_\nu$ exist, Theorem~\ref{id:general-BC} realizes the mixed Euler product as a product of same-field Rankin--Selberg functions. The residue of $K_{\Pi,\pi}^S$ at $1$ counts common cuspidal constituents of these lifts, with multiplicity. Every tensor-product factor is required, including when $F\subset E$.
\end{enumerate}
Higher-rank rational probes are included in the first construction: their base changes can contain several cuspidal constituents, and the residue records the multiplicity of the target. The fourth construction supplies a further mechanism over different fields, conditional on its stated base changes. The secondary terms contain more detailed expressions than these leading overlaps. The examples below make two limits explicit. First, local divisibility alone need not imply a global constituent: the symmetric-square example in Remark~\ref{id:divisibility-warning} has a local trivial root at every good prime but remains cuspidal. Second, zero statistics cannot distinguish $\Pi$ from a Galois conjugate with the same standard $L$-function. Thus the probe questions must be formulated for the grouped local polynomials, with their unavoidable symmetries, rather than for a recovery of labeled places.

\subsection{Pair, cubic, and quartic correlations}
We work out the pair, cubic, and quartic correlations explicitly. For a pair, Proposition~\ref{pa:finite} identifies the regular coefficients after subtracting the double pole of the pair kernel, and Corollary~\ref{cor:pair} places them in the zero covariance. Its first arithmetic correction is of order $(\log T)^{-2}$. The cubic constant in \eqref{hi:cubic-mixed} includes terms involving $b_\Pi(p)^2\overline{b_\pi(p^2)}$, so it joins degree-one data to degree-two residue data and second local powers. The quartic calculation subtracts all three pairings; \eqref{hi:quartic-mixed} exhibits both a mixed square and the product of two absolute squares. The full cubic and quartic constants, including every higher prime-power term, converge absolutely. They enter the centered connected statistics at orders $(\log T)^{-3}$ and $(\log T)^{-4}$, respectively. To display the additional local information, the lowest balanced-degree contributions for labels $(\Pi,\Pi,\pi)$ and $(\Pi,\Pi,\pi,\pi)$ are, respectively,
\begin{align}
&\frac{2(\log p)^3}{p^2}\Re\left(b_\Pi(p)^2\overline{b_\pi(p^2)}
+2b_\Pi(p)b_\pi(p)\overline{b_\Pi(p^2)}\right),\label{in:cubic-local}\\
&-\frac{(\log p)^4}{p^2}\left(2\Re\bigl(b_\Pi(p)^2\overline{b_\pi(p)^2}\bigr)
+4|b_\Pi(p)|^2|b_\pi(p)|^2\right).\label{in:quartic-local}
\end{align}
They are portions of the full absolutely convergent constants, not replacements for the higher-power terms. The distinct algebraic expressions in \eqref{in:cubic-local}--\eqref{in:quartic-local} explain why one should examine more than the leading pair overlap. Their injectivity as a complete set of functorial invariants remains an inverse problem.

More generally, put $\lambda=c\log T$ and center each zero statistic by its exact archimedean and pole term. For separable tests with $\operatorname{supp}\widehat f_j\subset[-\sigma_j,\sigma_j]$, we work under
\begin{equation}c\sum_j\sigma_j<1.\label{in:support}\end{equation}
Smooth height averaging then separates distinct prime-product frequencies. Theorem~\ref{cor:higher} gives, for every fixed $r\ge3$,
\begin{equation}\operatorname{cum}_T(X_{C_1},\ldots,X_{C_r})=\frac{(-1)^r\Gamma_{\boldsymbol C,\boldsymbol0}}{\lambda^r}\prod_j\widehat f_j(0)+O(\lambda^{-r-1}),\label{in:cumulant}\end{equation}
and every further finite inverse-logarithmic order, for arbitrary fixed $E,F$. This connected statement requires no mixed prime number theorem. Raw and distinct-index correlations follow from the partition reconstruction and sieving in Theorem~\ref{cor:raw-theorem} and \eqref{cor:sieve} when every required pair PNT is available. A mixed PNT by itself does not supply the self pair PNT over an arbitrary field. The support condition is an explicit sufficient box; no enlargement to a general joint Fourier-support region is claimed.

\subsection{Earlier work}
The history of this approach begins with Montgomery's pair correlation \cite{Montgomery}, Hejhal's triple correlation \cite{Hejhal}, and Rudnick--Sarnak's higher correlations \cite{RS}. Superpositions and their functorial interpretations were developed by Liu and Ye \cite{LiuYe2002,LiuYe2008,LiuYe2009}. The PNT of Liu and Ye \cite{LiuYe2007}, and subsequent work over number fields \cite{GillespieJi}, supplied coefficient orthogonality in their settings. Lower-order arithmetic in PNTs, Hypothesis~H, and zero correlations was studied by Gillespie and Ye \cite{GY2014,GY2016}; Euler--Stieltjes constants over Galois fields were studied by Od\v{z}ak and Smajlovi\'c \cite{OS2016}. Gillespie's work \cite{Gillespie2013,Gillespie2014,Gillespie2017} relates factorization, cross-field Rankin--Selberg constructions, zero statistics, and induction. The present contribution develops the local and secondary arithmetic forced by absolute norms for arbitrary fields and connects it to the precise detection statements above.

\subsection{Hypothesis H and the prime number theorem}
Two analytic inputs make these computations possible. First, we reprove Hypothesis~H for fixed representations over $E$ and also establish the grouped estimate
\begin{equation}\sum_p\sum_{k\ge2}\frac{(\log p)^2|b_\Pi(p^k)|^2}{p^{k(1-\eta)}}<\infty,\qquad 0<\eta<\frac1{d_Em^2+1}.\label{in:H}\end{equation}
This controls first local powers at residue degree greater than one as well as the usual higher prime-ideal powers. The proof uses a power sieve inspired by Jiang \cite{Jiang}, with ordinary Dirichlet characters pulled back by the norm, and applies also to products of positive Rankin--Selberg series. It supplies convergence with a positive margin for every higher balanced arithmetic block. Applying the positive-coefficient sieve directly to prime-ideal powers separately proves the usual formulation of Hypothesis~H; we do not deduce it by reversing the grouped coefficient sum. The distinction matters: at a place with $Nv=p^f$, $f>1$, the first local trace already occurs at a higher rational-prime power. A hypothesis restricted to local powers of exponent at least two does not cover it. Our proof avoids extending prescribed-order characters from a ray-class kernel. Appendix~\ref{app:jiang} gives a concrete obstruction to that extension in the cited version of Jiang's argument; it does not question the rational-field fixed-data argument. Second, the theorem of Harcos--Thorner \cite{HT} gives a PNT with an arbitrarily large logarithmic saving for each genuine same-field Rankin--Selberg factor, without a self-contragredience assumption. Any nonzero norm-twist main terms must be retained; the unshifted pair expansion is stated when such terms are absent. Thus the pair input is proved in the realizations stated above, rather than imposed there as an unproved conjecture. General unrelated-field products still require a realization or another analytic argument. All estimates allow dependence on the fixed fields and representations.

\subsection{Conventions}
Throughout, $n$ denotes the rank of the probe and $r$ the order of a correlation. For a representation $C$ over $E_C$, write $d_C=[E_C:\mathbb Q]$ and $\Delta_C=d_Cm_C$. Superscripts $S$ omit finitely many rational primes; their local factors are restored in the secondary constants. A tilde denotes the contragredient. Standard Rankin--Selberg properties are as in \cite{JPSS,JS,JSII}, and the local bounds are recalled in \cite[\S3.3]{HT}, following \cite{LRS,MS}.

\section{Rational norms and the mixed Euler product}\label{id:sec-coefficients}

Let $E$ and $F$ be number fields of degrees $d_E$ and $d_F$ over $\Q$.
Let $\Pi$ and $\pi$ be fixed unitary cuspidal automorphic representations of
$\mathrm{GL}_m(\mathbb A_E)$ and $\mathrm{GL}_{n}(\mathbb A_F)$, respectively.
Choose a finite set $S$ of rational primes containing the primes ramified in
either field and the rational primes below every ramified finite component
of either representation. A superscript $S$ means that all places above
$S$ are omitted. The set $S$ may be enlarged when additional representations
are introduced.

For $p\notin S$ and $v\mid p$ in $E$, write
\begin{equation}\label{id:local-data}
 Nv=p^{f_v},\qquad
 A_{\Pi,v}=\operatorname{diag}(\alpha_{v,1},\ldots,\alpha_{v,m}),\qquad
 a_\Pi(v^r)=\sum_{i=1}^{m}\alpha_{v,i}^{\,r}.
\end{equation}
The matrix notation denotes a semisimple conjugacy class; a particular
ordering of its eigenvalues is irrelevant. Define
\begin{equation}\label{id:grouped-coefficients}
 b_\Pi(p^k)=
 \sum_{\substack{v\mid p\\ f_v\mid k}}
 f_v\,a_\Pi(v^{k/f_v}).
\end{equation}
The analogous notation for $\pi$ will use places $w$ of $F$, residue degrees
$g_w$, and eigenvalues $\beta_{w,j}$. Direct expansion of the Euler product
gives
\begin{align}
 -\frac{(L_E^S)'}{L_E^S}(s,\Pi)
 &=\sum_{v\nmid S}\sum_{r\ge1}
   \frac{(\log Nv)a_\Pi(v^r)}{(Nv)^{rs}}\notag\\
 &=\sum_{p\notin S}\sum_{r\ge1}\sum_{v\mid p}
   \frac{f_v(\log p)a_\Pi(v^r)}{p^{f_vrs}}
 =\sum_{p\notin S}\sum_{k\ge1}
   \frac{(\log p)b_\Pi(p^k)}{p^{ks}}.
 \label{id:logarithmic-coefficients}
\end{align}
Thus the coefficient at a rational integer $p^k$ collects all prime-ideal
powers of that absolute norm. This grouping is forced by the exponential
frequencies in the explicit formula.

\subsection{A local polynomial without a global induction hypothesis}
At a good rational prime put
\begin{equation}\label{id:formal-polynomial}
 P_{\Pi,p}(X)=\prod_{v\mid p}\prod_{i=1}^{m}
                  (1-\alpha_{v,i}X^{f_v}).
\end{equation}
Since $\sum_{v\mid p}f_v=d_E$, this polynomial has degree $d_Em$.
Expanding each logarithm yields the formal identity
\begin{align}
 -\log P_{\Pi,p}(X)
 &=\sum_{v\mid p}\sum_{r\ge1}
               \frac{a_\Pi(v^r)X^{f_vr}}{r}\notag\\
 &=\sum_{k\ge1}\frac{b_\Pi(p^k)}{k}X^k.
 \label{id:polynomial-traces}
\end{align}
Equivalently, factor $P_{\Pi,p}(X)=\prod_{a=1}^{d_Em}(1-\gamma_{p,a}X)$
over $\mathbb C$. For each $v,i$, the associated $f_v$ numbers $\gamma_{p,a}$
are all the $f_v$th roots of $\alpha_{v,i}$. Their power sums are
\begin{equation}\label{id:formal-power-sums}
 \sum_{a=1}^{d_Em}\gamma_{p,a}^{\,k}=b_\Pi(p^k).
\end{equation}
The local multisets $\{\gamma_{p,a}\}$ are available whether or not they are
the Satake parameters of a global automorphic representation over $\Q$.
No such global representation is assumed in this paper.

\begin{remark}\label{id:quadratic-example}
Suppose $E/\Q$ is quadratic. At a split prime $p=v_1v_2$,
\[
 b_\Pi(p^k)=a_\Pi(v_1^k)+a_\Pi(v_2^k),\qquad
 P_{\Pi,p}(X)=\det(1-A_{\Pi,v_1}X)\det(1-A_{\Pi,v_2}X).
\]
At an inert prime there is one place $v$ of norm $p^2$, and
\[
 b_\Pi(p^{2r+1})=0,\qquad b_\Pi(p^{2r})=2a_\Pi(v^r),\qquad
 P_{\Pi,p}(X)=\det(1-A_{\Pi,v}X^2).
\]
In particular, the first local trace at an inert place occurs at the
\emph{rational} square $p^2$. This explains why a hypothesis involving only
prime-ideal powers of exponent at least two does not immediately control
every higher rational-prime power.
\end{remark}

\subsection{The mixed logarithmic series}
Define, initially for $\Re s>1$,
\begin{equation}\label{id:mixed-K}
 K_{\Pi,\pi}^S(s)=
 \sum_{p\notin S}\sum_{k\ge1}
 \frac{(\log p)b_\Pi(p^k)\overline{b_\pi(p^k)}}{p^{ks}}.
\end{equation}
We will justify this domain of absolute convergence below. The conjugation
agrees with passage to the contragredient because the representations are
unitary.

\begin{proposition}[Exact local factor]\label{id:local-factor-prop}
For $v\mid p$ in $E$ and $w\mid p$ in $F$, put
\[
 f=f_v,\qquad g=g_w,\qquad h_{v,w}=\gcd(f,g),\qquad
 L_{v,w}=\operatorname{lcm}(f,g).
\]
Then the Euler factors
\begin{equation}\label{id:mixed-local-factor}
 \mathscr L_{\Pi,\pi,p}(s)=
 \prod_{v\mid p}\prod_{w\mid p}
 \prod_{i=1}^{m}\prod_{j=1}^{n}
 \left(1-
 \alpha_{v,i}^{\,L_{v,w}/f_v}
 \overline{\beta_{w,j}^{\,L_{v,w}/g_w}}
 p^{-L_{v,w}s}\right)^{-h_{v,w}}
\end{equation}
satisfy
\begin{equation}\label{id:mixed-log-derivative}
 -\frac{d}{ds}\log\mathscr L_{\Pi,\pi}^S(s)
 =K_{\Pi,\pi}^S(s),\qquad
 \mathscr L_{\Pi,\pi}^S(s)=\prod_{p\notin S}\mathscr L_{\Pi,\pi,p}(s).
\end{equation}
The logarithm of the Euler product converges absolutely in $\Re s>1$.
\end{proposition}

\begin{proof}
Fix $v,w$ and abbreviate $h=h_{v,w}$, $L=L_{v,w}$. The logarithmic
derivative of their contribution to \eqref{id:mixed-local-factor} equals
\begin{align*}
 hL(\log p)\sum_{r\ge1}
 \frac{a_\Pi(v^{rL/f_v})
             \overline{a_\pi(w^{rL/g_w})}}{p^{rLs}}
 =f_vg_w(\log p)\sum_{\substack{k\ge1\\ f_v\mid k,\ g_w\mid k}}
 \frac{a_\Pi(v^{k/f_v})\overline{a_\pi(w^{k/g_w})}}{p^{ks}}.
\end{align*}
Here we used $hL=f_vg_w$, and then set $k=rL$. Summing over $v,w$ gives
\eqref{id:mixed-log-derivative} formally.

To prove convergence, weighted Cauchy--Schwarz gives
\begin{equation}\label{id:diagonal-majorant}
 |b_\Pi(p^k)|^2
 \le d_E\sum_{\substack{v\mid p\\ f_v\mid k}}
           f_v|a_\Pi(v^{k/f_v})|^2.
\end{equation}
Consequently, for $\sigma>1$,
\begin{align}
 \sum_{p\notin S}\sum_{k\ge1}
       \frac{|b_\Pi(p^k)|^2}{kp^{k\sigma}}
 &\le d_E\sum_{v\nmid S}\sum_{r\ge1}
             \frac{|a_\Pi(v^r)|^2}{r(Nv)^{r\sigma}}\notag\\
 &=d_E\log L_E^S(\sigma,\Pi\times\widetilde\Pi)<\infty.
 \label{id:diagonal-convergence}
\end{align}
The last equality uses unramified Rankin--Selberg Euler factors and their
nonnegative logarithmic coefficients. Applying Cauchy--Schwarz once more,
now to the sums over $(p,k)$, proves
\[
 \sum_{p\notin S}\sum_{k\ge1}
 \frac{|b_\Pi(p^k)\overline{b_\pi(p^k)}|}{kp^{k\sigma}}<\infty.
\]
This is the absolute convergence of $\log\mathscr L_{\Pi,\pi}^S$.
Logarithmic differentiation is justified locally uniformly in $\Re s>1$:
the extra factor $k\log p=\log(p^k)$ is absorbed by replacing $\sigma$ with
any fixed number between $1$ and $\sigma$ in the preceding bound.
\end{proof}

The construction is closely related to Gillespie's cross-field Euler
products \cite{Gillespie2014}. In the Galois setting, that work develops
both downstairs and upstairs constructions and connects their analytic
properties with mixed zero correlations and automorphic induction.
Formula \eqref{id:mixed-local-factor} allows varying residue degrees and
therefore does not require either field to be Galois. Absolute convergence
alone gives neither continuation across $\Re s=1$ nor a prime number
theorem for \eqref{id:mixed-K}.

\section{Rankin--Selberg realizations and probe detection}\label{id:sec-realizations}\label{id:section}

The identities in this section identify situations in which the mixed
Euler product has the analytic properties needed later. Each assumption
of automorphic base change is explicit. These assumptions concern the
probe or the field factors specified below, not automorphic induction of
$\Pi$ to $\Q$.

\subsection{A rational-field probe with known base change}
Suppose $F=\Q$. By a weak base change of $\pi$ to $E$ we mean an automorphic
representation $\tau$ whose local Satake class at every good $v\mid p$
has eigenvalues
\begin{equation}\label{id:BC-local-relation}
 \beta_{p,1}^{\,f_v},\ldots,\beta_{p,n}^{\,f_v}.
\end{equation}
For the arguments below we require $\tau$ to be a unitary isobaric sum of
cuspidal representations. The existence of such a lift is an input;
solvable base change supplies important settings in which it is known
\cite{AC}. The terminology does not impose a base-change theorem for an
arbitrary extension.

\begin{theorem}[Rational-field probe identity]\label{id:Q-probe}
Assume that $\tau=\operatorname{BC}_{E/\Q}(\pi)$ exists as a unitary isobaric
representation with \eqref{id:BC-local-relation}. Then
\begin{equation}\label{id:Q-probe-identity}
 \mathscr L_{\Pi,\pi}^S(s)=L_E^S(s,\Pi\times\widetilde\tau),\qquad
 K_{\Pi,\pi}^S(s)=-\frac{(L_E^S)'}{L_E^S}(s,\Pi\times\widetilde\tau).
\end{equation}
The pole order of the first function at $s=1$, equivalently the residue
of the second function there, is the multiplicity of $\Pi$ among the
cuspidal constituents of $\tau$.
\end{theorem}

\begin{proof}
For $F=\Q$ we have $g_w=1$, $L_{v,w}=f_v$, and $h_{v,w}=1$.
Thus \eqref{id:mixed-local-factor} becomes
\[
 \mathscr L_{\Pi,\pi,p}(s)=
 \prod_{v\mid p}\prod_{i=1}^{m}\prod_{j=1}^{n}
 (1-\alpha_{v,i}\overline{\beta_{p,j}^{\,f_v}}(Nv)^{-s})^{-1}.
\]
These are exactly the Rankin--Selberg factors of $\Pi\times\widetilde\tau$.
For completeness, one can also check the logarithmic coefficients directly:
\begin{align*}
 -\frac{(L_E^S)'}{L_E^S}(s,\Pi\times\widetilde\tau)
 &=\sum_{p\notin S}\sum_{v\mid p}\sum_{r\ge1}
 \frac{f_v(\log p)a_\Pi(v^r)\overline{a_\pi(p^{f_vr})}}{p^{f_vrs}}\\
 &=\sum_{p\notin S}\sum_{k\ge1}
 \frac{(\log p)b_\Pi(p^k)\overline{b_\pi(p^k)}}{p^{ks}}.
\end{align*}
Write $\tau=\boxplus_a\tau_a$, listing constituents with multiplicity.
Rankin--Selberg theory gives a simple pole at $s=1$ for
$L_E(s,\Pi\times\widetilde{\tau_a})$ precisely when $\tau_a\simeq\Pi$;
otherwise the factor is holomorphic and nonzero there. Finite omitted
Euler factors are nonzero at $1$. Multiplying the constituent factors
therefore proves the assertion about pole order.
\end{proof}

\begin{corollary}[An exact Satake relationship detected by a probe]
\label{id:Satake-detection}
Under the hypotheses of Theorem~\ref{id:Q-probe}, let
$\rho_{\Pi,\pi}=\operatorname*{Res}_{s=1}K_{\Pi,\pi}^S(s)$.
Then $\rho_{\Pi,\pi}$ is a nonnegative integer, and
$\rho_{\Pi,\pi}>0$ if and only if $\Pi$ is a constituent of $\tau$.
If $m=n$, this is equivalent to $\tau\simeq\Pi$, and it implies
\begin{equation}\label{id:Satake-detected-equality}
 \{\alpha_{v,1},\ldots,\alpha_{v,m}\}
 =\{\beta_{p,1}^{\,f_v},\ldots,\beta_{p,m}^{\,f_v}\}
 \qquad(v\mid p,\ p\notin S).
\end{equation}
Conversely, equality of these local multisets outside a finite set, with
the stated automorphic lift available, implies $\tau\simeq\Pi$ by strong
multiplicity one.
\end{corollary}

The equal-rank conclusion follows because a constituent of degree $m$
exhausts an isobaric representation of total degree $m$. In this case
$\rho_{\Pi,\pi}=1$. When the ranks differ, the pole criterion supplies a
global constituent relationship rather than merely an inclusion of
local eigenvalue multisets.

\begin{remark}[Local divisibility is weaker]\label{id:divisibility-warning}
Pointwise divisibility of Euler polynomials is not sufficient for a global
constituent assertion. Let $\varpi$ be a non-dihedral cuspidal
$\mathrm{GL}_2/\Q$ representation with trivial central character, so that
its symmetric-square lift is cuspidal \cite{GJ}. If the good-prime
parameters of $\varpi$ are $\alpha_p,\alpha_p^{-1}$, those of
$\operatorname{Sym}^2\varpi$ are
$\alpha_p^2,1,\alpha_p^{-2}$. The polynomial $1-X$ divides every good-prime
Euler polynomial of the symmetric-square lift, but the trivial
$\mathrm{GL}_1$ representation is not an isobaric constituent of that
cuspidal $\mathrm{GL}_3$ representation. The pole in
Theorem~\ref{id:Q-probe} supplies the global compatibility that local
divisibility by itself lacks.
\end{remark}

\subsection{Two representations over the same Galois field}
Suppose now that $E/\Q$ is Galois, with group $G$, and that $\Pi,\Sigma$
are unitary cuspidal representations over $E$, of possibly different
ranks. Define the conjugate $\Sigma^g$ by transporting local components,
so that
$a_{\Sigma^g}(v^r)=a_\Sigma((g^{-1}v)^r)$ at good finite places.

\begin{theorem}[Galois conjugate identity]\label{id:Galois-identity-thm}
Without assuming automorphic induction or any nontrivial base change,
\begin{equation}\label{id:Galois-identity}
 \mathscr L_{\Pi,\Sigma}^S(s)=
 \prod_{g\in G}L_E^S(s,\Pi\times\widetilde{\Sigma^g}).
\end{equation}
Consequently,
\begin{equation}\label{id:Galois-residue}
 \operatorname*{Res}_{s=1}K_{\Pi,\Sigma}^S(s)
 =\#\{g\in G:\Pi\simeq\Sigma^g\}.
\end{equation}
In particular,
\begin{equation}\label{id:stabilizer}
 \operatorname*{Res}_{s=1}K_{\Pi,\Pi}^S(s)
 =\#\operatorname{Stab}_G(\Pi).
\end{equation}
\end{theorem}

\begin{proof}
At $p\notin S$, all places above $p$ have the same residue degree $f$.
The stabilizer in $G$ of any such place has order $f$, because $p$ is
unramified. Hence, for any fixed $v\mid p$, every $w\mid p$ occurs exactly
$f$ times in the list $g^{-1}v$, $g\in G$. It follows that
\[
 \sum_{g\in G}\overline{a_{\Sigma^g}(v^r)}
 =f\sum_{w\mid p}\overline{a_\Sigma(w^r)}.
\]
The logarithmic derivative of the right side of
\eqref{id:Galois-identity}, at the frequency $p^{fr}$, therefore has
coefficient
\[
 (\log p)f\sum_{v\mid p}a_\Pi(v^r)
       \sum_{g\in G}\overline{a_{\Sigma^g}(v^r)}
 =f^2(\log p)
   \sum_{v,w\mid p}a_\Pi(v^r)\overline{a_\Sigma(w^r)}.
\]
This is $(\log p)b_\Pi(p^{fr})\overline{b_\Sigma(p^{fr})}$.
Both logarithmic derivatives vanish at exponents not divisible by $f$.
Thus their coefficients agree. The Euler products tend to $1$ as
$\Re s\to+\infty$, which fixes the constant after integration and proves
\eqref{id:Galois-identity}. The pole criterion applied to each factor
gives \eqref{id:Galois-residue} and \eqref{id:stabilizer}.
\end{proof}

\begin{remark}\label{id:Galois-self-shifts}
There are no nonzero norm-twist coincidences in the self-case of
\eqref{id:Galois-identity}. Indeed, if
$\Pi^g\simeq\Pi\otimes|\det|^{it}$ and $g$ has order $a$, iteration gives
$\Pi\simeq\Pi\otimes|\det|^{iat}$. Comparing central characters on the
positive real norm subgroup gives $t=0$. The self-PNT below consequently
has the single main term
$\#\operatorname{Stab}_G(\Pi)\,x$.
\end{remark}

\begin{remark}[A split quadratic norm collision]\label{id:split-collision}
Let $E/\Q$ be quadratic and $p=v_1v_2$ split. The rational mixed
coefficient is
\begin{align*}
 b_\Pi(p^k)\overline{b_\Sigma(p^k)}
 ={}&a_\Pi(v_1^k)\overline{a_\Sigma(v_1^k)}
    +a_\Pi(v_2^k)\overline{a_\Sigma(v_2^k)}\\
 &+a_\Pi(v_1^k)\overline{a_\Sigma(v_2^k)}
    +a_\Pi(v_2^k)\overline{a_\Sigma(v_1^k)}.
\end{align*}
The first line comes from $L_E(s,\Pi\times\widetilde\Sigma)$; the second
comes from the nontrivial Galois conjugate. At an inert prime,
$b_\Pi(p^{2r})\overline{b_\Sigma(p^{2r})}
=4a_\Pi(v^r)\overline{a_\Sigma(v^r)}$.
There are two conjugate Rankin--Selberg factors, each with logarithmic
weight $\log Nv=2\log p$, again giving the required factor $4\log p$.
\end{remark}

\subsection{The tensor product of the two ground fields}
There is a uniform formulation of both preceding identities. The finite
etale $\Q$-algebra $E\otimes_\Q F$ is a product of number fields:
\begin{equation}\label{id:etale-decomposition}
 E\otimes_\Q F\simeq\prod_{\nu=1}^{u}K_\nu.
\end{equation}
Each factor has two specified embeddings, one from $E$ and one from $F$.
All base changes in the next theorem refer to these embeddings.

\begin{theorem}[General base-change realization]\label{id:general-BC}
Assume that for every factor in \eqref{id:etale-decomposition} there are
unitary isobaric weak base changes
\[
 \tau_\nu=\operatorname{BC}_{K_\nu/E}(\Pi),\qquad
 \upsilon_\nu=\operatorname{BC}_{K_\nu/F}(\pi).
\]
Then, after enlarging $S$ if necessary,
\begin{equation}\label{id:general-BC-identity}
 \mathscr L_{\Pi,\pi}^S(s)=
 \prod_{\nu=1}^{u}
 L_{K_\nu}^S(s,\tau_\nu\times\widetilde{\upsilon_\nu}).
\end{equation}
\end{theorem}

\begin{proof}
Fix a good $p$ and a pair $v\mid p$, $w\mid p$ of residue degrees $f,g$.
The extensions $E_v/\Q_p$ and $F_w/\Q_p$ are unramified. Their tensor
product consists of $h=\gcd(f,g)$ copies of the unramified extension of
$\Q_p$ of degree $L=\operatorname{lcm}(f,g)$. This follows, for example,
by the corresponding finite-field tensor-product decomposition, followed
by unramified lifting. The local degree over $E_v$ is $L/f$ and the local
degree over $F_w$ is $L/g$ in each copy. The respective base-change
Satake eigenvalues are therefore
\[
 \alpha_{v,i}^{\,L/f},\qquad \beta_{w,j}^{\,L/g}.
\]
The Rankin--Selberg factor in one copy is
\[
 \prod_{i,j}
 \left(1-\alpha_{v,i}^{\,L/f}
       \overline{\beta_{w,j}^{\,L/g}}p^{-Ls}\right)^{-1}.
\]
Multiplying over the $h$ copies and then over $v,w$ gives exactly
\eqref{id:mixed-local-factor}. This proves the identity prime by prime.
\end{proof}

When $F\subset E$, one factor of \eqref{id:etale-decomposition} is $E$,
but the other factors generally cannot be discarded. The single
Rankin--Selberg function obtained by base-changing $\pi$ along the
specified inclusion $F\subset E$ accounts only for compatible places.
The remaining factors account for other pairs of places with equal
absolute norms. When $E=F$ is Galois, the decomposition has one copy of
$E$ for every $g\in G$ and reduces to
Theorem~\ref{id:Galois-identity-thm}.

\subsection{The precise prime number theorem supplied by these identities}
For representations $\rho,\sigma$ over the same number field $K$, put
\[
 \mathcal M_{\rho,\sigma}(x)=
 \begin{cases}
 x^{1+it}/(1+it),&\sigma\simeq\rho\otimes|\det|^{it}
                         \text{ for some }t\in\mathbb R,\\
 0,&\text{otherwise}.
 \end{cases}
\]
Here the Rankin--Selberg function under consideration is
$L_K(s,\rho\times\widetilde\sigma)$; with this convention a twist by
$|\det|^{it}$ in $\sigma$ produces a pole at $s=1+it$.
Harcos--Thorner \cite[Theorem 2.1]{HT}, summed over the narrow ideal
classes at modulus $(1)$, gives
\begin{equation}\label{id:HT-PNT}
 \sum_{N\mathfrak a\le x}
 \Lambda_{\rho\times\widetilde\sigma}(\mathfrak a)
 =\mathcal M_{\rho,\sigma}(x)
   +O_{K,\rho,\sigma,A}\bigl(x(\log x)^{-A}\bigr)
 \qquad(A>0).
\end{equation}
This is a same-field statement, without a self-contragredience
assumption. Its implied constant is ineffective. It does not assert a
general subexponential remainder. The narrow class group at modulus
$(1)$ need not be trivial; summing over its classes cancels its
nontrivial characters and gives the unrestricted ideal sum above.

\begin{corollary}[Mixed PNT in the realized cases]\label{id:mixed-PNT}
Suppose one of the Rankin--Selberg realizations proved above is available.
List every cuspidal constituent pair $(\rho_j,\sigma_j)$ in the
corresponding product, retaining multiplicity and its ground field.
Then, for every $A>0$,
\begin{equation}\label{id:mixed-PNT-formula}
 \sum_{\substack{p^k\le x\\p\notin S}}
  (\log p)b_\Pi(p^k)\overline{b_\pi(p^k)}
 =\sum_j\mathcal M_{\rho_j,\sigma_j}(x)
       +O_{\Pi,\pi,E,F,A}\bigl(x(\log x)^{-A}\bigr).
\end{equation}
If no constituent pair is related by a nonzero norm twist, the main term
is $\rho_{\Pi,\pi}x$, where
\begin{equation}\label{id:residue-general}
 \rho_{\Pi,\pi}
 =\sum_j\boldsymbol 1_{\rho_j\simeq\sigma_j}
 =\operatorname*{Res}_{s=1}K_{\Pi,\pi}^S(s).
\end{equation}
\end{corollary}

\begin{proof}
Logarithmic differentiation of the relevant realization identifies its
coefficients, grouped by rational norm, with the left side of
\eqref{id:mixed-PNT-formula}. Apply \eqref{id:HT-PNT} to each of the
finitely many constituent pairs and add the resulting formulas.
\end{proof}

\begin{remark}[Finite Euler factors]\label{id:finite-factors}
The omission of $S$ does not change the residues above. To check this
quantitatively, the standard local root bounds give
$|\alpha_{v,i}|\le(Nv)^{\theta_m}$ with
$\theta_m\le\tfrac12-(m^2+1)^{-1}$, and the corresponding bound for
$\pi$. Hence the contribution of a fixed rational prime to
\eqref{id:mixed-K} is holomorphic for
$\Re s>\theta_m+\theta_{n}<1$. The genuine Rankin--Selberg factors at
the omitted places likewise have local exponents strictly below $1$;
see \cite[\S3.3]{HT}. For some fixed $\vartheta<1$, the coefficients
removed from any of the PNT sums have total size $O(x^\vartheta)$,
after increasing $\vartheta$ slightly if necessary. This is absorbed
by every error in \eqref{id:mixed-PNT-formula}. The omitted factors
must nevertheless be restored when computing finite parts and secondary
constants.
\end{remark}

\section{A power sieve for rational norms}\label{sec:norm-sieve}

The coefficients entering zero correlations are indexed by rational integers. Consequently, an estimate restricted to ideals of the form $\mathfrak b^k$ is not quite the estimate needed here. We prove an estimate for the larger set of ideals whose absolute norms are $k$th powers. This distinction is essential at places of residue degree greater than one.

The argument uses the power-sieve idea underlying Jiang's approach to Hypothesis~H \cite{Jiang}; compare also the rational power sieve in \cite{Munshi}. We apply ordinary Dirichlet characters to absolute norms. This avoids any extension of prescribed-order characters from the kernel of a ray class group. Appendix~\ref{app:jiang} explains why that extension cannot be taken for granted. All fields and automorphic representations are fixed throughout this section.

\subsection{A positive Dirichlet series over the rational integers}

Let $E_1,\ldots,E_J$ be number fields, with $d_j=[E_j:\Q]$, and let $\Pi_j$ be a unitary cuspidal automorphic representation of $\mathrm{GL}_{m_j}(\mathbb A_{E_j})$. Repetitions among the fields or representations are allowed. Let $S$ be a finite set of rational primes containing all primes ramified in any $E_j$ and all rational primes below a finite place where some $\Pi_j$ is ramified. We use the superscript $S$ to omit every local factor above these rational primes. Put
\begin{equation}\label{ns:positive-series}
\mathcal D(s)=\prod_{j=1}^{J}L_{E_j}^{S}(s,\Pi_j\times\widetilde\Pi_j)
             =\sum_{N\ge1}\frac{A(N)}{N^s},
\qquad D=\sum_{j=1}^{J}d_jm_j^2.
\end{equation}
The coefficients $A(N)$ are nonnegative. More explicitly,
\begin{equation}\label{ns:coeff-convolution}
A(N)=\sum_{\substack{\mathfrak a_j\subseteq\mathcal O_{E_j}\ (1\le j\le J)\\
                    \prod_jN\mathfrak a_j=N\\
                    (N\mathfrak a_j,S)=1\ (1\le j\le J)}}
       \prod_{j=1}^{J}\lambda_{\Pi_j\times\widetilde\Pi_j}(\mathfrak a_j).
\end{equation}
At an unramified place, the positivity also follows directly from
\begin{equation}\label{ns:positive-local-exp}
\sum_{h\ge0}\lambda_{\Pi_j\times\widetilde\Pi_j}(v^h)z^h
 =\exp\left(\sum_{h\ge1}\frac{|a_{\Pi_j}(v^h)|^2}{h}z^h\right).
\end{equation}
Here unitarity identifies the multiset of contragredient Satake roots with the complex conjugate multiset. Standard Rankin--Selberg positivity, including the general local formulation, is recalled in \cite[\S3.4]{HT}.

Each diagonal Rankin--Selberg factor has a simple pole at $1$. The omitted local factors are finite and nonzero there, so $\mathcal D(s)$ has a pole of order $J$. In particular, for $0<u\le1$,
\begin{equation}\label{ns:pole-bound}
\mathcal D(1+u)\ll u^{-J}.
\end{equation}
All constants in this section may depend on the fixed fields, representations, and $S$.

\begin{theorem}[Norm-power estimate]\label{ns:norm-power}
Let $A(N)$ and $D$ be as in \eqref{ns:positive-series}. Fix an integer $k\ge2$, a nonnegative function $V\in C_c^\infty((0,\infty))$, and $\varepsilon>0$. Then
\begin{equation}\label{ns:norm-power-bound}
\sum_{\substack{N\ge1\\N\in\mathbb Z_{\ge1}^{\,k}}}A(N)V(N/X)
\ll X^{1-1/(D+1)+\varepsilon}\qquad(X\ge2).
\end{equation}
The implied constant is allowed to depend on $k$, $V$, $\varepsilon$, and all fixed automorphic and field data. No uniform assertion in discriminants or automorphic conductors is made.
\end{theorem}

\begin{proof}
We give the character construction, the sieve inequality, and the analytic estimate separately.

\emph{Step 1: rational primes and norm characters.}
Let $M$ be the compositum of the normal closures of $E_1,\ldots,E_J$ and the cyclotomic field $\Q(\zeta_{2k})$. For a parameter $P$, let
\begin{equation}\label{ns:sieving-primes}
\mathcal P(P)=\{q:P<q\le2P,\ q\notin S,\ q\text{ splits completely in }M\},
\qquad R=\#\mathcal P(P).
\end{equation}
The fixed-field Chebotarev theorem gives
\begin{equation}\label{ns:prime-count}
R\asymp_{M}\frac{P}{\log P}
\end{equation}
for sufficiently large $P$. In particular, each $q\in\mathcal P(P)$ satisfies $q\equiv1\pmod{2k}$.

Choose a Dirichlet character $\chi_q$ modulo $q$ of exact order $k$, extended by zero on multiples of $q$. It is primitive because $q$ is prime and $k\ge2$. Moreover it is even: if a generator of $\mathbb F_q^\times$ is sent to a primitive $k$th root of unity, its $(q-1)/2$th power is sent to $1$, since $k\mid(q-1)/2$. For each $j$, pull back its associated Hecke character by the idelic norm:
\begin{equation}\label{ns:norm-character}
\theta_{q,j}=\chi_q\circ N_{E_j/\Q},\qquad
\theta_{q,j}(\mathfrak a)=\chi_q(N\mathfrak a)
\end{equation}
on ideals coprime to $q$. This construction is globally well-defined by the norm map on idele class groups. Its infinite component is trivial because $\chi_q$ is even. At every place $v\mid q$, complete splitting identifies $(E_j)_v$ with $\Q_q$, and the local norm is the identity. Thus $\theta_{q,j}$ has finite conductor $q\mathcal O_{E_j}$.

For distinct $q,r\in\mathcal P(P)$, put
\begin{equation}\label{ns:cross-character}
\psi_{q,r}=\chi_q\overline{\chi_r},\qquad
\Theta_{q,r,j}=\theta_{q,j}\overline{\theta_{r,j}}.
\end{equation}
The character $\Theta_{q,r,j}$ has conductor $qr\mathcal O_{E_j}$ and trivial infinite component. Its ramification at every place above $q$ and $r$ will rule out all possible Rankin--Selberg poles.

\emph{Step 2: the elementary sieve inequality.}
Choose fixed $0<a<b$ with $\operatorname{supp}V\subseteq[a,b]$, and put $w(N)=A(N)V(N/X)\ge0$. If $N$ is a $k$th power, then
\begin{equation}\label{ns:power-detector}
\sum_{q\in\mathcal P(P)}\chi_q(N)
=R-\#\{q\in\mathcal P(P):q\mid N\}.
\end{equation}
For $N\le bX$, the number subtracted is at most $\log(bX)/\log P$. We will take $P=X^{1/(D+1)}$. In view of \eqref{ns:prime-count}, for all sufficiently large $X$ this number is at most $R/2$. It follows that
\begin{align}
\sum_{N\in\mathbb Z_{\ge1}^{\,k}}w(N)
&\le\frac4{R^2}\sum_{N\ge1}w(N)
               \left|\sum_{q\in\mathcal P(P)}\chi_q(N)\right|^2\notag\\
&\le\frac4R\sum_{N\ge1}w(N)
 +4\max_{\substack{q,r\in\mathcal P(P)\\q\ne r}}
       \left|\sum_{N\ge1}w(N)\psi_{q,r}(N)\right|.
\label{ns:sieve-inequality}
\end{align}
The second line follows by expanding the square. Each of the $R$ diagonal terms is at most $\sum_Nw(N)$; there are at most $R^2$ off-diagonal terms. This calculation uses positivity of $w$, but makes no assertion that the off-diagonal sums themselves are positive.

For the diagonal sum, take $u=1/\log(bX)$, with $X$ sufficiently large. Rankin's trick and \eqref{ns:pole-bound} give
\begin{align}
\sum_Nw(N)&\le\|V\|_\infty\sum_{N\le bX}A(N)
 \le\|V\|_\infty(bX)^{1+u}\mathcal D(1+u)\notag\\
&\ll X(\log X)^J.
\label{ns:diagonal-bound}
\end{align}
In particular, even for a single diagonal factor the direct estimate furnished by this argument is $O(X\log X)$. It is not legitimate to replace the pole bound $\mathcal D(1+u)\ll u^{-1}$ by a bounded quantity as $u\downarrow0$. The logarithm has no effect on the final exponent.

\emph{Step 3: the off-diagonal Dirichlet series.}
For $\Re s>1$, complete multiplicativity of $\psi_{q,r}$ and \eqref{ns:coeff-convolution} give the exact identity
\begin{align}
\mathcal D_{q,r}(s)
&:=\sum_{N\ge1}\frac{A(N)\psi_{q,r}(N)}{N^s}\notag\\
&=\prod_{j=1}^{J}L_{E_j}^{S}
 \bigl(s,\Pi_j\times(\widetilde\Pi_j\otimes\Theta_{q,r,j})\bigr).
\label{ns:twisted-series}
\end{align}
At places not above $q,r$ or $S$, the equality is immediate from the unramified Euler factors. At a place above $q$ or $r$, the Rankin--Selberg parameter before twisting is unramified, whereas the twisting character is ramified. Its inertia-fixed subspace is therefore zero and the twisted local L-factor is $1$. This agrees exactly with the convention $\psi_{q,r}(N)=0$ when $(N,qr)>1$; no missing local correction occurs at the moving primes.

A pole of a factor on the right of \eqref{ns:twisted-series} could occur only if
\begin{equation}\label{ns:pole-exclusion}
\widetilde\Pi_j\otimes\Theta_{q,r,j}
\simeq\widetilde\Pi_j\otimes|\det|^{it}
\end{equation}
for a real $t$. At a place above $q$, the left side is a ramified twist of an unramified representation, whereas the right side is unramified. Their local Langlands parameters have different inertia actions. Hence \eqref{ns:pole-exclusion} is impossible. The completed Rankin--Selberg factors consequently have no poles, including shifted poles. The finite L-functions in \eqref{ns:twisted-series} are holomorphic in the strip used below.

The conductor exponent at each place above $q$ or $r$ is $m_j^2$. Thus, with fixed archimedean data,
\begin{equation}\label{ns:twisted-conductor}
C\bigl(it,\Pi_j\times(\widetilde\Pi_j\otimes\Theta_{q,r,j})\bigr)
\ll (qr)^{d_jm_j^2}(1+|t|)^{d_jm_j^2}.
\end{equation}
Constants in this bound incorporate the fixed discriminant and the fixed ramified data. Rankin--Selberg convexity, for example \cite[Lemma~3.2]{HT}, gives for fixed $0<\sigma<1$ and every $\epsilon_1>0$
\begin{equation}\label{ns:convexity}
|\mathcal D_{q,r}(\sigma+it)|
\ll_{\sigma,\epsilon_1}
P^{D(1-\sigma)+\epsilon_1}(1+|t|)^{B_{\sigma,\epsilon_1}}.
\end{equation}
Here $B_{\sigma,\epsilon_1}$ is a fixed exponent. Removing the factors above $S$ multiplies the complete finite L-functions by finitely many local inverse polynomials. On a fixed vertical line these polynomials are bounded uniformly in $q,r$, since their roots have fixed absolute values and the character values have absolute value one. This justifies applying convexity to the incomplete functions in \eqref{ns:convexity}.

Let $\widehat V(s)=\int_0^\infty V(x)x^{s-1}\,dx$. It is entire and decays faster than any power of $|\Im s|$ in each fixed vertical strip. Mellin inversion gives
\begin{equation}\label{ns:mellin}
\sum_{N\ge1}A(N)\psi_{q,r}(N)V(N/X)
=\frac1{2\pi i}\int_{(2)}\mathcal D_{q,r}(s)X^s\widehat V(s)\,ds.
\end{equation}
Shift the line to $\Re s=\epsilon_2$, where $\epsilon_2>0$ is small. No pole is crossed. Equation~\eqref{ns:convexity}, together with the decay of $\widehat V$, yields
\begin{equation}\label{ns:off-diagonal}
\sum_NA(N)\psi_{q,r}(N)V(N/X)
\ll P^{D+\epsilon_1}X^{\epsilon_2}.
\end{equation}

\emph{Step 4: choice of the sieve parameter.}
Combining \eqref{ns:prime-count}, \eqref{ns:sieve-inequality}, \eqref{ns:diagonal-bound}, and \eqref{ns:off-diagonal}, we obtain
\begin{equation}\label{ns:pre-optimization}
\sum_{N\in\mathbb Z_{\ge1}^{\,k}}A(N)V(N/X)
\ll\frac{X(\log X)^J\log P}{P}+P^{D+\epsilon_1}X^{\epsilon_2}.
\end{equation}
With $P=X^{1/(D+1)}$, the two principal powers of $X$ are both $X^{D/(D+1)}$. Choose $\epsilon_1,\epsilon_2$ sufficiently small in terms of the prescribed $\varepsilon$ and absorb the logarithms. This proves \eqref{ns:norm-power-bound} for sufficiently large $X$. Enlarging the fixed implied constant covers the remaining bounded range.
\end{proof}

\begin{remark}\label{ns:fixed-data-scope}
The proof uses only fixed-field Chebotarev, ordinary Dirichlet characters, standard Rankin--Selberg analytic properties, and convexity. It neither assumes automorphic induction nor uses a mixed-field prime number theorem. Because the Chebotarev field $M$ depends on the original fields and on $k$, the proof makes no claim of the uniformity in field discriminants or automorphic conductors asserted in stronger formulations of Jiang's theorem.
\end{remark}

\subsection{The rationally grouped form of Hypothesis H}

Now take a single representation $\Pi$ of $\mathrm{GL}_m(\mathbb A_E)$, and write $d=[E:\Q]$. At $p\notin S$, define
\begin{equation}\label{ns:grouped-coeff}
b_\Pi(p^r)=\sum_{\substack{v\mid p\\f_v\mid r}}
 f_v\,a_\Pi(v^{r/f_v}),\qquad Nv=p^{f_v}.
\end{equation}
Let $A_\Pi(N)$ be the coefficients of $L_E^S(s,\Pi\times\widetilde\Pi)$ grouped by absolute norm, as in \eqref{ns:positive-series} with $J=1$.

\begin{lemma}\label{ns:coefficient-domination}
For every $p\notin S$ and integer $r\ge1$,
\begin{equation}\label{ns:coefficient-bound}
|b_\Pi(p^r)|^2\le dr\,A_\Pi(p^r).
\end{equation}
\end{lemma}
\begin{proof}
The coefficient of $z^h$ on the right of \eqref{ns:positive-local-exp} is a sum of nonnegative terms, one of which is $|a_\Pi(v^h)|^2/h$. Hence
\begin{equation}\label{ns:local-domination}
|a_\Pi(v^h)|^2\le h\lambda_{\Pi\times\widetilde\Pi}(v^h).
\end{equation}
Weighted Cauchy--Schwarz and $\sum_{v\mid p}f_v=d$ at a good rational prime give
\begin{align*}
|b_\Pi(p^r)|^2
&\le\left(\sum_{\substack{v\mid p\\f_v\mid r}}f_v\right)
      \left(\sum_{\substack{v\mid p\\f_v\mid r}}
                   f_v|a_\Pi(v^{r/f_v})|^2\right)\\
&\le d\sum_{\substack{v\mid p\\f_v\mid r}}
     f_v\frac r{f_v}\lambda_{\Pi\times\widetilde\Pi}(v^{r/f_v})\\
&=dr\sum_{\substack{v\mid p\\f_v\mid r}}
      \lambda_{\Pi\times\widetilde\Pi}(v^{r/f_v})
 \le dr A_\Pi(p^r).
\end{align*}
For the last inequality, the displayed prime-ideal powers are distinct ideals of norm $p^r$, and all occur among the nonnegative summands defining $A_\Pi(p^r)$.
\end{proof}

\begin{theorem}[Grouped Hypothesis H with a positive margin]\label{ns:collapsed-h}
Put
\begin{equation}\label{ns:delta-definition}
\delta_\Pi=\frac1{dm^2+1}.
\end{equation}
For every $0<\eta<\delta_\Pi$,
\begin{equation}\label{ns:collapsed-h-bound}
\sum_p\sum_{r\ge2}
\frac{(\log p)^2|b_\Pi(p^r)|^2}{p^{r(1-\eta)}}<\infty.
\end{equation}
Here the coefficients at the finitely many bad rational primes may be defined from the standard local L-factors and included in the sum. For every fixed nonnegative integer $B$, the same conclusion holds after inserting a factor $r^B(1+\log p)^B$, provided $\eta<\delta_\Pi$ remains strict.
\end{theorem}

\begin{proof}
First fix $r\ge2$. Choose $V\in C_c^\infty((0,\infty))$ with $V\ge1$ on $[1,2]$. Theorem~\ref{ns:norm-power}, applied with $k=r$, and Lemma~\ref{ns:coefficient-domination} give
\begin{equation}\label{ns:dyadic-fixed-r}
\sum_{\substack{p\notin S\\X\le p^r\le2X}}|b_\Pi(p^r)|^2
\le dr\sum_{N\in\mathbb Z_{\ge1}^{\,r}}A_\Pi(N)V(N/X)
\ll_r X^{1-\delta_\Pi+\epsilon}.
\end{equation}
Choose $0<\epsilon<\delta_\Pi-\eta$. On this dyadic interval, $(\log p)^2\ll_r(\log X)^2$, so
\begin{equation}\label{ns:dyadic-weighted}
\sum_{\substack{p\notin S\\X\le p^r\le2X}}
\frac{(\log p)^2|b_\Pi(p^r)|^2}{p^{r(1-\eta)}}
\ll_r X^{\eta-\delta_\Pi+\epsilon}(\log X)^2.
\end{equation}
Summing over $X=2^u$ proves convergence for each fixed $r$. This step does not require constants uniform in $r$.

It remains to sum over all sufficiently large $r$. The Luo--Rudnick--Sarnak local root bound, with its general local formulation recalled in \cite[\S3.3]{HT}, gives
\begin{equation}\label{ns:root-bound}
|\alpha_{\Pi,v,i}|\le(Nv)^{1/2-a_m},\qquad a_m=\frac1{m^2+1}.
\end{equation}
At ramified places the list of local roots may be padded by zeros. Since $\sum_{v\mid p}f_v\le d$, equation~\eqref{ns:grouped-coeff} and the same formula using local Euler roots at bad primes imply
\begin{equation}\label{ns:pointwise-grouped}
|b_\Pi(p^r)|\le dm\,p^{r(1/2-a_m)}.
\end{equation}
As $\eta<\delta_\Pi\le a_m$, the number $2a_m-\eta$ is positive. Choose an integer $R\ge2$ with $R(2a_m-\eta)>1$. Then
\begin{align}
\sum_p\sum_{r\ge R}
\frac{(\log p)^2|b_\Pi(p^r)|^2}{p^{r(1-\eta)}}
&\le d^2m^2\sum_p(\log p)^2
       \sum_{r\ge R}p^{-r(2a_m-\eta)}\notag\\
&=d^2m^2\sum_p(\log p)^2
\frac{p^{-R(2a_m-\eta)}}{1-p^{-(2a_m-\eta)}}<\infty.
\label{ns:large-r-tail}
\end{align}
There are only finitely many exponents $2\le r<R$, already treated by \eqref{ns:dyadic-weighted}; for these exponents the finitely many bad primes contribute finite sums. This proves \eqref{ns:collapsed-h-bound}.

For the final assertion choose $\eta'$ with $\eta<\eta'<\delta_\Pi$. Since $p\ge2$, every fixed polynomial in $r$ and $\log p$ is bounded by a constant times $p^{r(\eta'-\eta)}$. Apply the result with $\eta'$.
\end{proof}

\begin{remark}\label{ns:residue-degree-warning}
The term $a_\Pi(v)$ at a place with residue degree $f_v\ge2$ occurs in $b_\Pi(p^{f_v})$. Its local exponent is one, whereas its rational exponent is at least two. It is not covered merely by asserting a prime-ideal hypothesis for exponents $h\ge2$. The norm-power sieve covers it directly because $Nv=p^{f_v}$ is a rational $f_v$th power. Conversely, the same argument also implies the usual prime-ideal Hypothesis H: for each fixed $h\ge2$, the terms $v^h$ have perfect-power norms, and \eqref{ns:local-domination} followed by Theorem~\ref{ns:norm-power} proves their summability. Thus the present estimate supplies both the customary hypothesis and the stronger grouped estimate used here.
\end{remark}

\begin{corollary}[Mixed grouped estimate]\label{ns:mixed-h}
Let $\Pi/\mathrm{GL}_m(E)$ and $\pi/\mathrm{GL}_{n}(F)$ be fixed unitary cuspidal representations, and put
\[
\delta_\Pi=\frac1{[E:\Q]m^2+1},\qquad
\delta_\pi=\frac1{[F:\Q]n^2+1}.
\]
For every $0<\eta<(\delta_\Pi+\delta_\pi)/2$,
\begin{equation}\label{ns:mixed-h-bound}
\sum_p\sum_{r\ge2}
\frac{(\log p)^2|b_\Pi(p^r)\overline{b_\pi(p^r)}|}
     {p^{r(1-\eta)}}<\infty.
\end{equation}
Arbitrary fixed polynomial factors in $r$ and $\log p$ may again be inserted.
\end{corollary}
\begin{proof}
Choose $0<\eta_\Pi<\delta_\Pi$ and $0<\eta_\pi<\delta_\pi$ such that
$\eta_\Pi+\eta_\pi>2\eta$, and set $\eta_*=(\eta_\Pi+\eta_\pi)/2$. The Cauchy--Schwarz inequality gives
\begin{align*}
&\sum_{p,r\ge2}\frac{(\log p)^2|b_\Pi(p^r)b_\pi(p^r)|}{p^{r(1-\eta_*)}}\\
&\qquad\le
\left(\sum_{p,r\ge2}\frac{(\log p)^2|b_\Pi(p^r)|^2}{p^{r(1-\eta_\Pi)}}\right)^{1/2}
\left(\sum_{p,r\ge2}\frac{(\log p)^2|b_\pi(p^r)|^2}{p^{r(1-\eta_\pi)}}\right)^{1/2}<\infty.
\end{align*}
The result at $\eta<\eta_*$ follows by comparison. The strict gap also absorbs any fixed polynomial factors, exactly as at the end of the proof of Theorem~\ref{ns:collapsed-h}.
\end{proof}

\section{The pair kernel and its secondary arithmetic terms}\label{pa:section}
The arithmetic series in this section use all rational primes. At a bad prime the logarithmic coefficients are defined from the actual local standard $L$-factor. The good-prime identities of Section~\ref{id:section} therefore acquire explicit finite Euler corrections when they are used here. This convention is essential for secondary terms.

For two labels $C,D$, standing for fixed unitary cuspidal representations over possibly different fields, put
\begin{align}
 K_{C,D}(s)&=\sum_{p,k\ge1}\frac{(\log p)b_C(p^k)\overline{b_D(p^k)}}{p^{ks}},\label{pa:K}\\
 A_{C,D}(s)&=\sum_{p,k\ge1}\frac{(\log p)^2b_C(p^k)\overline{b_D(p^k)}}{p^{ks}},\label{pa:A}\\
 H_{C,D}(s)&=\sum_{p,k\ge2}\frac{(k-1)(\log p)^2b_C(p^k)\overline{b_D(p^k)}}{p^{ks}}.\label{pa:H}
\end{align}
Absolute convergence in $\Re s>1$ permits differentiation, and gives
\begin{equation} A_{C,D}(s)=-K_{C,D}'(s)-H_{C,D}(s).\label{pa:identity}\end{equation}
The factor $k-1$ in \eqref{pa:H} cannot be omitted: differentiating $p^{-ks}$ produces $k\log p$, whereas the explicit formula supplies only $\log p$ in each of its two factors.

Let $\delta_C=([E_C:\mathbb Q]m_C^2+1)^{-1}$. The norm-power estimates proved above imply that $H_{C,D}$ is holomorphic in
\begin{equation}\Re s>1-\tfrac12(\delta_C+\delta_D).\label{pa:Hstrip}\end{equation}
For clarity, holomorphy on this open half-plane means local uniform convergence: on any compact subset one first chooses a strictly larger convergence margin. The resulting exponential factor absorbs $k-1$ and any fixed power of $k\log p$. This also justifies every derivative of $H_{C,D}$ used below.

\begin{proposition}[Arithmetic finite part]\label{pa:finite}
Suppose $K_{C,D}$ has a meromorphic continuation to a neighborhood of $1$, with expansion
\begin{equation}K_{C,D}(1+z)=\frac{\rho_{C,D}}z+\sum_{j\ge0}\kappa_{C,D,j}z^j.\label{pa:Kexp}\end{equation}
Then
\begin{equation}A_{C,D}(1+z)=\frac{\rho_{C,D}}{z^2}+\sum_{j\ge0}C_{C,D,j}z^j,\qquad
C_{C,D,j}=-(j+1)\kappa_{C,D,j+1}-\frac{H_{C,D}^{(j)}(1)}{j!}.\label{pa:Cexp}\end{equation}
In particular,
\begin{equation}C_{C,D,0}=-\kappa_{C,D,1}-H_{C,D}(1).\label{pa:C0}\end{equation}
There is no simple-pole term in $A_{C,D}$.
\end{proposition}
\begin{proof}
Insert \eqref{pa:Kexp} into \eqref{pa:identity}. The derivative of $\rho/z$ is $-\rho/z^2$, the derivative of $\kappa_0$ is zero, and the coefficient of $z^j$ in the derivative of the regular part is $(j+1)\kappa_{j+1}$. Taylor expansion of the holomorphic function $H$ gives \eqref{pa:Cexp}.
\end{proof}

We next give a test-function version that specifies exactly how these constants enter a pair statistic. It is useful independently of the zero-correlation application.

\begin{theorem}[Smoothed pair arithmetic]\label{pa:smooth}
Assume, for every $B>0$, that
\begin{equation}S_{C,D}(x):=\sum_{p^k\le x}(\log p)b_C(p^k)\overline{b_D(p^k)}
=\rho_{C,D}x+O_{C,D,B}\bigl(x(\log(2x))^{-B}\bigr).\label{pa:PNT}\end{equation}
Let $\phi$ be smooth and compactly supported on $[0,\infty)$, with right derivatives at zero, and put
\begin{equation}\mathcal D_{C,D}(\lambda;\phi)=\frac1{\lambda^2}\sum_{p,k\ge1}
\frac{(\log p)^2b_C(p^k)\overline{b_D(p^k)}}{p^k}\,
\phi\left(\frac{k\log p}{\lambda}\right).\label{pa:D}\end{equation}
For every fixed integer $J\ge0$, as $\lambda\to\infty$,
\begin{equation}\mathcal D_{C,D}(\lambda;\phi)
=\rho_{C,D}\int_0^\infty u\phi(u)\,du+
\sum_{j=0}^{J}\frac{(-1)^j C_{C,D,j}\phi^{(j)}(0)}{\lambda^{j+2}}
+O_{C,D,\phi,J}(\lambda^{-J-3}).\label{pa:Dexp}\end{equation}
When \eqref{pa:Kexp} is available, the constants are exactly those in \eqref{pa:Cexp}. Under \eqref{pa:PNT} alone they are equivalently defined by the convergent moment integrals in the proof. If $\phi$ is constant near zero, then for every $N>0$,
\begin{equation}\mathcal D_{C,D}(\lambda;\phi)
=\rho_{C,D}\int_0^\infty u\phi(u)\,du+
\frac{C_{C,D,0}\phi(0)}{\lambda^2}+O_{C,D,\phi,N}(\lambda^{-N}).\label{pa:plateau}\end{equation}
\end{theorem}
\begin{proof}
Write $R(x)=S_{C,D}(x)-\rho_{C,D}x$ for $x\ge1$ and set
\begin{equation}M_j=\int_1^\infty R(x)x^{-2}(\log x)^j\,dx\qquad(j\ge0).\label{pa:Mj}\end{equation}
These integrals converge absolutely: in \eqref{pa:PNT} choose $B>j+1$ and put $y=\log x$. The resulting integrand is $O((1+y)^{j-B})\,dy$.

First replace the second $\log p$ in \eqref{pa:D} by $\log(p^k)=k\log p$. Before division by $\lambda^2$, this gives
\begin{equation}\int_{1^-}^\infty x^{-1}\log x\,
\phi\left(\frac{\log x}{\lambda}\right)dS_{C,D}(x).\label{pa:Stieltjes}\end{equation}
The contribution of $\rho x$ is exactly $\rho\lambda^2\int_0^\infty u\phi(u)du$. For the remainder, integration by parts is legitimate because the weight vanishes at $x=1$ and at the upper endpoint of its support. It gives
\begin{equation}-\int_1^\infty\frac{R(x)}{x^2}
\left\{(1-\log x)\phi\left(\frac{\log x}{\lambda}\right)
+\frac{\log x}{\lambda}\phi'\left(\frac{\log x}{\lambda}\right)\right\}dx.\label{pa:remainder}\end{equation}
Taylor expansion at zero shows that the coefficient of $\lambda^{-j}$ in \eqref{pa:remainder} is
\begin{equation}-\frac{\phi^{(j)}(0)}{j!}\bigl((j+1)M_j-M_{j+1}\bigr).\label{pa:Rcoefficient}\end{equation}
The error after $j=J$ is $O(\lambda^{-J-1})$, because the Taylor remainders are bounded by constants times the corresponding powers of $(\log x)/\lambda$, and all the integrals \eqref{pa:Mj} are finite.

For $\Re z>0$, partial summation also gives
\begin{equation}K_{C,D}(1+z)=\frac{\rho_{C,D}}z+\rho_{C,D}
+(1+z)\int_1^\infty R(x)x^{-2-z}dx.\label{pa:Kmoments}\end{equation}
Consequently
\begin{equation}\kappa_{C,D,j+1}
=\frac{(-1)^jM_j}{j!}+\frac{(-1)^{j+1}M_{j+1}}{(j+1)!}.
\label{pa:kappamoments}\end{equation}
The identity holds for Taylor coefficients when the continuation exists; otherwise the right side defines the required boundary derivatives. Combining \eqref{pa:Rcoefficient} and \eqref{pa:kappamoments} gives the contribution $(-1)^j[-(j+1)\kappa_{C,D,j+1}]\phi^{(j)}(0)$.

Finally subtract the excess higher powers,
\begin{equation}\sum_{p,k\ge2}\frac{(k-1)(\log p)^2b_C(p^k)\overline{b_D(p^k)}}{p^k}
\phi\left(\frac{k\log p}{\lambda}\right).\label{pa:Hweight}\end{equation}
Absolute convergence with a positive margin permits Taylor expansion term by term. Its coefficient of $\lambda^{-j}$ is
$(-1)^jH_{C,D}^{(j)}(1)\phi^{(j)}(0)/j!$, with error $O(\lambda^{-J-1})$. Divide the resulting identity by $\lambda^2$ to obtain \eqref{pa:Dexp}. If $\phi$ is constant near zero, all its positive-order derivatives there vanish. Choose $J$ arbitrarily large to obtain \eqref{pa:plateau}.
\end{proof}

\begin{remark}[Shifted coincidences]
The hypothesis \eqref{pa:PNT} is written without nonzero norm-twist main terms. In the Rankin--Selberg situations of Section~\ref{id:section}, it follows from Harcos--Thorner when no such coincidences occur. In general the finitely many terms $\rho_t x^{1+it}/(1+it)$, $t\ne0$, must be kept separately. Their contribution to \eqref{pa:Stieltjes}, after division by $\lambda^2$, is
\[
\rho_t\int_0^\infty u\phi(u)e^{it\lambda u}\,du.
\]
Repeated integration by parts gives its complete inverse-power expansion. Thus a nonzero shift contributes to the regular arithmetic terms at $s=1$; it does not add to $\rho_{C,D}$. Constants in this expansion are not uniform as $t\to0$. We state subsequent pair formulas under \eqref{pa:PNT} to avoid suppressing this distinction.
\end{remark}

\begin{remark}[Error terms and earlier work]
The all-logarithmic PNT error suffices for every finite expansion \eqref{pa:Dexp}. It does not, by itself, imply a subexponential remainder for arbitrary smooth test functions. Lower-order arithmetic constants in orthogonality and zero correlations were already studied by Gillespie and Ye \cite{GY2014,GY2016}; the purpose here is to retain the rational-norm grouping for different fields and to identify the constants relevant to probes.
\end{remark}

\section{Higher arithmetic blocks and mixed cumulants}\label{hi:section}

The pair kernel is exceptional because its first-prime terms need a prime
number theorem.  The arithmetic blocks of every fixed order at least three
are different: the balance of their frequencies forces a perfect-power norm,
and Theorem~\ref{ns:norm-power} gives absolute convergence.  We prove this
assertion before introducing cumulants, so that every rearrangement below is
justified.

Throughout this section, $\mathcal C=\{\Pi,\pi\}$ is a list of the two
representations; it is regarded as a list even when its entries coincide.
For $C\in\mathcal C$, write $d_C$ for the degree of its field and $m_C$ for
its rank.  Put
\begin{equation}\label{hi:parameters}
 \begin{gathered}
 D=d_Em^2+d_Fn^2,\qquad \delta_*=(D+1)^{-1},\qquad
 d_* =\max(d_E,d_F),\\
 a_* =\min\{(m^2+1)^{-1},(n^2+1)^{-1}\}.
 \end{gathered}
\end{equation}
The set $S$ contains all the bad rational primes.  Initially every prime sum
is over $p\notin S$.  At the end of the section we explain the finite
corrections.

\subsection{A positive coefficient majorant}

Let $A(N)$ be the nonnegative rational-norm coefficients of
\begin{equation}\label{hi:positive-product}
 L_E^S(s,\Pi\times\widetilde\Pi)
 L_F^S(s,\pi\times\widetilde\pi)
 =\sum_{N\geq1}\frac{A(N)}{N^s}.
\end{equation}
At a good rational prime $p$, define
\begin{equation}\label{hi:ck}
 c_k(p)=\sum_{C\in\mathcal C}
 \sum_{\substack{v\mid p\text{ in the field of }C\\f_v\mid k}}
 f_v\,|a_C(v^{k/f_v})|^2.
\end{equation}
The local logarithm of \eqref{hi:positive-product} gives the identity of
power series
\begin{equation}\label{hi:exp-majorant}
 \sum_{K\geq0}A(p^K)z^K
 =\exp\left(\sum_{k\geq1}\frac{c_k(p)}{k}z^k\right).
\end{equation}
In particular all the coefficients on the right are nonnegative.  Weighted
Cauchy--Schwarz, using $\sum_{v\mid p}f_v=d_C$ at good primes, gives
\begin{equation}\label{hi:b-majorant}
 |b_C(p^k)|^2
 \leq d_C\sum_{\substack{v\mid p\\f_v\mid k}}
 f_v|a_C(v^{k/f_v})|^2
 \leq d_Cc_k(p).
\end{equation}

\begin{lemma}\label{hi:monomial-majorant}
Let $u,v\geq1$, $r=u+v$, and let $C_1,\ldots,C_u,D_1,\ldots,D_v$
be entries of $\mathcal C$.  If
\[
 k_1+\cdots+k_u=l_1+\cdots+l_v=K,\qquad k_i,l_j\geq1,
\]
then
\begin{equation}\label{hi:balanced-majorant}
 \left|\prod_{i=1}^u b_{C_i}(p^{k_i})
       \prod_{j=1}^v\overline{b_{D_j}(p^{l_j})}\right|
 \leq d_*^{r/2}K^{r/2}\sqrt{u!v!}\,A(p^K).
\end{equation}
\end{lemma}

\begin{proof}
For the positive exponents let $\mu_h=\#\{i:k_i=h\}$.  Thus
$\sum_h\mu_h=u$ and $\sum_hh\mu_h=K$.  Selecting the corresponding term
from each exponential factor in \eqref{hi:exp-majorant} yields
\[
 A(p^K)\geq\prod_{h\geq1}
 \frac{c_h(p)^{\mu_h}}{h^{\mu_h}\mu_h!}.
\]
Consequently, by \eqref{hi:b-majorant},
\begin{align*}
 \prod_{i=1}^u|b_{C_i}(p^{k_i})|^2
 &\leq d_*^u\prod_hc_h(p)^{\mu_h}\\
 &\leq d_*^u A(p^K)\prod_hh^{\mu_h}\mu_h!
 \leq d_*^uK^u u!\,A(p^K).
\end{align*}
Here $h\leq K$ and $\prod_h\mu_h!\leq u!$.  The same calculation for the
negative exponents gives $d_*^vK^v v!A(p^K)$.  Taking the square root of
the product proves \eqref{hi:balanced-majorant}.
\end{proof}

\subsection{Convergence with a positive margin}

For $\epsilon\in\{+1,-1\}$ set
\[
 b_C^{+1}(p^k)=b_C(p^k),\qquad
 b_C^{-1}(p^k)=\overline{b_C(p^k)}.
\]
Fix $r\geq3$, labels $C_1,\ldots,C_r\in\mathcal C$, and signs
$\epsilon_1,\ldots,\epsilon_r\in\{+1,-1\}$.  If all signs are the same,
the balance condition below has no solutions and the block is zero.
Otherwise define
\begin{equation}\label{hi:block}
 \mathcal B_{\boldsymbol C,\boldsymbol\epsilon}(\boldsymbol z)
 =\sum_{p\notin S}(\log p)^r
 \sum_{\substack{k_1,\ldots,k_r\geq1\\
                  \epsilon_1k_1+\cdots+\epsilon_rk_r=0}}
 \frac{\prod_{j=1}^r b_{C_j}^{\epsilon_j}(p^{k_j})}
 {p^{\frac12\sum_jk_j+\sum_jk_jz_j}}.
\end{equation}
The conjugation in this formula applies only to the coefficients.  Each
$z_j$ is an independent complex variable, not a conjugated variable.

\begin{theorem}\label{hi:convergence}
Let $r\geq3$ be fixed and let
\begin{equation}\label{hi:eta-range}
 0<\eta<\min\{\delta_*,2a_*\}.
\end{equation}
For every fixed $B,q\geq0$, the series
\begin{equation}\label{hi:absolute-margin}
 \sum_{p\notin S}(\log p)^B
 \sum_{\substack{k_j\geq1\\\sum_j\epsilon_jk_j=0}}
 \frac{(1+\sum_jk_j)^q\prod_j|b_{C_j}(p^{k_j})|}
 {p^{(1-\eta)\sum_jk_j/2}}
\end{equation}
converges.  In particular \eqref{hi:block} is holomorphic on
$|\Re z_j|<\eta/2$ for $1\leq j\leq r$, and all its partial derivatives
are given by termwise differentiation there.
\end{theorem}

\begin{proof}
Let $u$ and $v$ be the number of positive and negative signs.  On a balanced
tuple, both sums of exponents equal an integer $K$, and $\sum_jk_j=2K$.
Since $u+v\geq3$, one has $K\geq\max(u,v)\geq2$.

First fix $K$.  There are finitely many relevant exponent tuples.  By
Lemma~\ref{hi:monomial-majorant}, it suffices to prove
\begin{equation}\label{hi:fixedK}
 \sum_{p\notin S}\frac{(\log p)^B A(p^K)}{p^{K(1-\eta)}}<\infty.
\end{equation}
Choose a nonnegative smooth function $V$ supported in $[1/2,5/2]$ and at
least one on $[1,2]$.  Theorem~\ref{ns:norm-power}, applied to
\eqref{hi:positive-product} with power exponent $K$, gives
\[
 \sum_{X\leq p^K\leq2X}A(p^K)
 \leq\sum_{N\text{ a }K\text{th power}}A(N)V(N/X)
 \ll_{K,\varepsilon}X^{1-\delta_*+\varepsilon}.
\]
The contribution to \eqref{hi:fixedK} from this interval is at most
\[
 \ll_{K,\varepsilon}(\log(2X))^B
 X^{\eta-\delta_*+\varepsilon}.
\]
Taking $0<\varepsilon<\delta_*-\eta$ and summing over dyadic $X$ proves
\eqref{hi:fixedK}.  Constants may depend on $K$ at this stage; no uniformity
in the power-sieve exponent is being asserted or used.

For large $K$ we use a separate estimate.  The local root bounds imply
\begin{equation}\label{hi:root-bound}
 |b_C(p^k)|\leq d_Cm_Cp^{k(1/2-a_*)}.
\end{equation}
If $M=\max_{C\in\mathcal C}d_Cm_C$, a balanced monomial, after division by
$p^{K(1-\eta)}$, is therefore at most
$M^r p^{-K(2a_*-\eta)}$.  The number of tuples with given $K$ equals
\[
 \binom{K-1}{u-1}\binom{K-1}{v-1}\ll_r K^{r-2},
\]
with the usual value zero when an exponent composition is impossible.
Choose $K_0$ so large that $K_0(2a_*-\eta)>1$.  The portion of
\eqref{hi:absolute-margin} with $K\geq K_0$ is bounded by
\[
 \ll_{r,q}M^r\sum_p(\log p)^B
 \sum_{K\geq K_0}K^{r-2+q}p^{-K(2a_*-\eta)}<\infty.
\]
The inner series is a geometric series with a polynomial factor, uniformly
for $p\geq2$; its first exponent is strictly larger than one.  Together with
the finitely many applications of \eqref{hi:fixedK}, this proves convergence.

On a closed substrip $|\Re z_j|\leq\eta'/2$ with $\eta'<\eta$, the absolute
value of each denominator in \eqref{hi:block} is at least
$p^{K(1-\eta')}$.  A partial derivative of total order $q$ only inserts a
factor bounded by $(2K)^q(\log p)^q$.  The already proved estimate, with a
slightly larger margin than $\eta'$, gives locally uniform convergence of
every differentiated series.  This proves the holomorphy assertion.
\end{proof}

\begin{remark}\label{hi:bad-primes}
At a fixed bad prime the standard local factors still give
$|b_C(p^k)|\leq M_Cp^{k(1/2-a_C)}$ with $a_C>0$; zero local roots may be
included.  Therefore all the exponent sums in this section converge
geometrically at that prime in a sufficiently small shift neighborhood.
The finitely many bad-prime corrections can be reinstated.  They change
the constants, but do not affect the convergence conclusions.
\end{remark}

\subsection{Local phase averages and connected constants}

The following phase model is a convenient way to organize the arithmetic;
it is not an assertion that the actual zeros form this random process.
For a fixed prime, define
\begin{equation}\label{hi:local-phase}
 X_{C,p}(\theta;z)
 =(\log p)\sum_{k\geq1}
 \frac{b_C(p^k)e^{ik\theta}
       +\overline{b_C(p^k)}e^{-ik\theta}}
 {p^{k(1/2+z)}},\qquad 0\leq\theta\leq2\pi.
\end{equation}
This series converges uniformly in $\theta$ on compact subsets of
$\Re z>-a_C$, by the local bound.  In particular it defines a real-valued
function when $z=0$, and its mean over $\theta$ is zero.

Write $\langle Y\rangle_p=(2\pi)^{-1}\int_0^{2\pi}Y(\theta)\,d\theta$.
For a finite index set $I$, let
\[
 M_{I,p}(\boldsymbol z)
 =\left\langle\prod_{j\in I}X_{C_j,p}(\theta;z_j)\right\rangle_p.
\]
Integrating $e^{i(\sum_j\epsilon_jk_j)\theta}$ gives one when the frequency
is zero and zero otherwise.  Thus the sum over $p$ of $M_{I,p}$ is exactly
the sum of the $2^{|I|}$ blocks \eqref{hi:block}.  The empty sign configurations
contribute zero.

Let $\mathfrak P_r$ denote the set of partitions of $\{1,\ldots,r\}$.
Define the local joint cumulant by the finite partition formula
\begin{equation}\label{hi:cumulant}
 \kappa_p(C_1,\ldots,C_r;\boldsymbol z)
 =\sum_{\mathcal P\in\mathfrak P_r}
 (-1)^{|\mathcal P|-1}(|\mathcal P|-1)!
 \prod_{I\in\mathcal P}M_{I,p}(\boldsymbol z).
\end{equation}
Every partition containing a singleton contributes zero.  In every
remaining term the frequency balance holds separately within each block,
and hence also for the whole tuple.  Its expansion is a subsum of balanced
monomials of order $r$.  Theorem~\ref{hi:convergence} applies to that expansion,
including its absolute values.

\begin{corollary}\label{hi:cumulant-convergence}
For every fixed $r\geq3$, the connected arithmetic constant
\begin{equation}\label{hi:global-cumulant}
 \mathcal K(C_1,\ldots,C_r;\boldsymbol z)
 =\sum_{p\notin S}\kappa_p(C_1,\ldots,C_r;\boldsymbol z)
\end{equation}
is absolutely convergent and holomorphic in a neighborhood of
$\boldsymbol z=0$.  All logarithmic moments obtained by partial
differentiation converge absolutely.
\end{corollary}

For completeness, take independent uniform phases $\theta_p$ only for
$p\leq Y$, and form the finite sum
$X_C^{(Y)}(z)=\sum_{p\leq Y,\ p\notin S}X_{C,p}(\theta_p;z)$.
Additivity of joint cumulants of independent vectors gives
\[
 \operatorname{cum}(X_{C_1}^{(Y)}(z_1),\ldots,X_{C_r}^{(Y)}(z_r))
 =\sum_{p\leq Y,\ p\notin S}\kappa_p(C_1,\ldots,C_r;\boldsymbol z).
\]
Corollary~\ref{hi:cumulant-convergence} permits $Y\to\infty$ on the right
for $r\geq3$.  We do not define an unregularized infinite random sum at
$z=0$: its variance can diverge.  Indeed, in the half-plane where the
prime sum converges, the covariance is
\begin{equation}\label{hi:covariance}
 \sum_{p\notin S}M_{\{1,2\},p}(z_1,z_2)
 =A_{C_1,C_2}^S(1+z_1+z_2)+A_{C_2,C_1}^S(1+z_1+z_2).
\end{equation}
This displays the pair singularity separately from the convergent higher
connected constants.

\subsection{The cubic calculation}

Because each local mean is zero, the third cumulant equals the third
moment.  Fix labels $C_1,C_2,C_3$ and put all shifts equal to zero.
Exactly two frequencies have one sign and the third has the opposite
sign.  For each $h\in\{1,2,3\}$ let $i<j$ be the two other indices.
The complete calculation is
\begin{equation}\label{hi:cubic}
 \kappa_p(C_1,C_2,C_3;\boldsymbol0)
 =2(\log p)^3\Re\sum_{h=1}^3\sum_{u,v\geq1}
 \frac{b_{C_i}(p^u)b_{C_j}(p^v)
       \overline{b_{C_h}(p^{u+v})}}{p^{u+v}}.
\end{equation}
The factor two includes the opposite sign configuration.  With labels
$\Pi,\Pi,\pi$, this becomes
\begin{equation}\label{hi:cubic-mixed}
 \begin{split}
 \kappa_p(\Pi,\Pi,\pi;\boldsymbol0)
 =2(\log p)^3\Re\sum_{u,v\geq1}\frac1{p^{u+v}}
 \bigl(&b_\Pi(p^u)b_\Pi(p^v)\overline{b_\pi(p^{u+v})}\\
 &+2b_\Pi(p^u)b_\pi(p^v)\overline{b_\Pi(p^{u+v})}\bigr).
 \end{split}
\end{equation}
In particular the $K=2$ contribution is
\begin{equation}\label{hi:cubic-first}
 \frac{2(\log p)^3}{p^2}\Re\left(
 b_\Pi(p)^2\overline{b_\pi(p^2)}
 +2b_\Pi(p)b_\pi(p)\overline{b_\Pi(p^2)}\right).
\end{equation}
The coefficients at $p^2$ include both second powers at degree-one places
and first powers at degree-two places.  Thus the cubic constant explicitly
mixes information from the two residue degrees.

\subsection{The quartic subtraction}

At order four, the connected constant is
\begin{equation}\label{hi:quartic-partitions}
 \kappa_p(1,2,3,4)=M_{1234,p}
 -M_{12,p}M_{34,p}-M_{13,p}M_{24,p}-M_{14,p}M_{23,p},
\end{equation}
where the labels and zero shifts are suppressed in this display only.
We compute its entire $K=2$ portion.  All four exponents then equal one.
Write $u_j=b_{C_j}(p)$ and $a=(\log p)/\sqrt p$.  The relevant first
harmonic is $a(u_je^{i\theta}+\overline{u_j}e^{-i\theta})$.
Its fourth moment is
\[
 a^4\sum_{\substack{I\subset\{1,2,3,4\}\\|I|=2}}
 \prod_{i\in I}u_i\prod_{j\notin I}\overline{u_j}.
\]
For each fixed two-element set $I$, the same monomial occurs twice among
the three covariance products in \eqref{hi:quartic-partitions}: there are
exactly two ways to pair its two positive frequencies with its two
negative frequencies.  Therefore the cumulant is the negative of the
displayed fourth moment.  Equivalently,
\begin{equation}\label{hi:quartic-first}
 -\frac{2(\log p)^4}{p^2}\Re\left(
 u_1u_2\overline{u_3u_4}
 +u_1u_3\overline{u_2u_4}
 +u_1u_4\overline{u_2u_3}\right).
\end{equation}
For labels $\Pi,\Pi,\pi,\pi$, this is
\begin{equation}\label{hi:quartic-mixed}
 -\frac{(\log p)^4}{p^2}
 \left(2\Re\bigl(b_\Pi(p)^2\overline{b_\pi(p)^2}\bigr)
       +4|b_\Pi(p)|^2|b_\pi(p)|^2\right).
\end{equation}
When all four labels coincide, it reduces to
$-6(\log p)^4|b_\Pi(p)|^4/p^2$.  The terms with $K\geq3$ must still be
included in the full quartic constant; their convergence follows from
Theorem~\ref{hi:convergence}.

\subsection{What these computations establish}

The results of this section establish the existence, shift holomorphy,
and explicit local formulas of higher arithmetic constants for arbitrary
fixed fields.  They use neither automorphic induction nor a mixed PNT.
They do not alone give the full raw correlation of distinct zeros.  The
explicit formula must still identify which arithmetic blocks occur and
control the terms with unequal prime-product frequencies.  In a sufficiently
small Fourier box, smooth averaging can perform this analytic step for
centered linear statistics; the resulting statement is proved separately
below.  A theorem for raw correlations of distinct zeros additionally
requires the archimedean contributions and the appropriate partition
sieving.  The phase model above records the exact balanced arithmetic
terms, and does not by itself supply those additional steps.

\section{Mixed zero statistics in a fixed Fourier box}\label{cor:section}
We now connect the convergent arithmetic blocks to zeros. The support condition below is deliberately explicit and sufficient for the argument. No enlargement of the Rudnick--Sarnak support region is claimed. Unlike a formal prime-phase analogy, the comparison in this section follows from the explicit formula and an elementary separation bound for distinct rational numbers.

\subsection{Normalization and the exact center}
For a label $C$ over $E_C$, put $\Delta_C=[E_C:\mathbb Q]m_C$; this is the degree of its completed standard $L$-function over $\mathbb Q$. Write
\begin{equation}\Lambda_C(s)=Q_C^{s/2}\prod_{\nu=1}^{\Delta_C}\Gamma_{\mathbb R}(s+\mu_{C,\nu})L(s,C),\qquad
\Gamma_{\mathbb R}(s)=\pi^{-s/2}\Gamma(s/2).\label{cor:completed}\end{equation}
Thus $Q_C=|D_{E_C}|^{m_C}N\mathfrak q_C$. We use the meromorphic completion, without multiplying by factors to cancel its possible poles. Its nontrivial zeros are denoted by $\rho_C$, counted with multiplicity, and
\begin{equation}\gamma_C=\frac{\rho_C-1/2}{i}.\label{cor:ordinate}\end{equation}
These ordinates need not be real. No generalized Riemann hypothesis is assumed.

Our Fourier convention is
\begin{equation}\widehat f(u)=\int_{\mathbb R}f(x)e^{-2\pi iux}\,dx.\label{cor:Fourier}\end{equation}
Take $\widehat f\in C_c^\infty(\mathbb R)$, with $f$ real on $\mathbb R$. Its inverse Fourier transform is entire, so evaluations at \eqref{cor:ordinate} are meaningful. Fix $c>0$, set $\lambda=c\log T$, and define
\begin{equation}Z_C(t;f)=\sum_{\rho_C}f\left(\frac{\lambda(\gamma_C-t)}{2\pi}\right).\label{cor:Z}\end{equation}
All sums over zeros in this section include multiplicities.

Let $g_t(u)=f(\lambda(u-t)/(2\pi))$, and define the exact archimedean and pole term
\begin{align}
 M_C(t;f)&=\sum_{\omega\text{ pole of }\Lambda_C}m_\omega\,
       g_t\left(\frac{\omega-1/2}{i}\right)\notag\\
 &\quad+\frac1{2\pi}\int_{\mathbb R}g_t(u)
 \left\{\log Q_C+
 2\Re\sum_{\nu=1}^{\Delta_C}
 \frac{\Gamma_{\mathbb R}'}{\Gamma_{\mathbb R}}
       \left(\frac12+iu+\mu_{C,\nu}\right)\right\}du.
 \label{cor:center}
\end{align}
Here $m_\omega$ denotes the order of the pole. The explicit formula, in the normalization of \cite[Proposition 2.1]{RS}, gives the exact identity
\begin{align}
 X_C(t;f)&:=Z_C(t;f)-M_C(t;f)\notag\\
 &=-\frac1\lambda\sum_{p,k\ge1}\frac{\log p}{p^{k/2}}
 \left\{b_C(p^k)\widehat f\left(\frac{k\log p}{\lambda}\right)e^{-itk\log p}
 +\overline{b_C(p^k)}\widehat f\left(-\frac{k\log p}{\lambda}\right)e^{itk\log p}\right\}.
 \label{cor:explicit}
\end{align}
The prime sum is finite. This formula includes the actual factors at bad primes. For real $f$, the two displayed terms are conjugate, so $X_C(t;f)$ is real.

Fix $w\in C_c^\infty((1,2))$ with $w\ge0$ and $\int w=1$. Write
\begin{equation}\mathbb E_T U=\frac1T\int_{\mathbb R}w(t/T)U(t)dt.\label{cor:height-average}\end{equation}
Joint cumulants with respect to this probability measure are denoted by $\operatorname{cum}_T$. Their definition is the same partition formula as \eqref{hi:cumulant}.

\subsection{An exact diagonal model up to a rapidly decreasing error}
Let $C_1,\ldots,C_r$ be fixed labels, each equal to $\Pi$ or $\pi$, and suppose
\begin{equation}\operatorname{supp}\widehat f_j\subseteq[-\sigma_j,\sigma_j],\qquad
\sigma_j>0,\qquad c\sum_{j=1}^r\sigma_j<1.\label{cor:support}\end{equation}
For independent phases $\theta_p$ uniform on $[0,2\pi]$, put
\begin{align}
Y_{j,p}^{(\lambda)}(\theta_p)
 =\sum_{k\ge1}\frac{\log p}{p^{k/2}}
 \biggl\{&b_{C_j}(p^k)\widehat f_j\left(\frac{k\log p}{\lambda}\right)e^{ik\theta_p}\notag\\
 &+\overline{b_{C_j}(p^k)}\widehat f_j\left(-\frac{k\log p}{\lambda}\right)e^{-ik\theta_p}\biggr\}.
 \label{cor:weighted-phase}
\end{align}
Only finitely many primes and exponents occur. In particular the random variables below are finite sums, not unregularized infinite series.

\begin{theorem}[Small-box comparison]\label{cor:comparison}
Under \eqref{cor:support}, for every $M>0$,
\begin{equation}\mathbb E_T\prod_{j=1}^rX_{C_j}(t;f_j)
=\frac{(-1)^r}{\lambda^r}\,
\mathbb E_{\boldsymbol\theta}\prod_{j=1}^r
 \left(\sum_pY_{j,p}^{(\lambda)}(\theta_p)\right)+O(T^{-M}).\label{cor:moment-comparison}\end{equation}
The implied constant depends on the fixed data, tests, support margin, $w,c,r,M$. Consequently
\begin{equation}\operatorname{cum}_T(X_{C_1}(t;f_1),\ldots,X_{C_r}(t;f_r))
=\frac{(-1)^r}{\lambda^r}\sum_p
\operatorname{cum}_{\theta_p}(Y_{1,p}^{(\lambda)},\ldots,Y_{r,p}^{(\lambda)})
+O(T^{-M}).\label{cor:cumulant-comparison}\end{equation}
Neither conclusion requires a mixed prime number theorem.
\end{theorem}
\begin{proof}
Expand the finite product in \eqref{cor:explicit}. A summand has frequency
\begin{equation}\Omega=\sum_{j=1}^r\epsilon_jk_j\log p_j=\log(a/b),\qquad \epsilon_j\in\{1,-1\},\label{cor:frequency}\end{equation}
where $a$ and $b$ are positive integers. Each nonzero Fourier weight implies
$p_j^{k_j}\le e^{\lambda\sigma_j}=T^{c\sigma_j}$. Hence
\begin{equation}\max(a,b)\le T^\vartheta,\qquad \vartheta=c\sum_j\sigma_j<1.\label{cor:ab-bound}\end{equation}
If $a\ne b$, then
\begin{equation}|\Omega|=|\log(a/b)|\ge\frac1{\max(a,b)}\ge T^{-\vartheta}.\label{cor:spacing}\end{equation}
For example, if $a>b$, integrate $1/u$ on $[b,a]$ to get $\log(a/b)\ge(a-b)/a\ge1/a$.

The height average of the frequency is
\begin{equation}\frac1T\int w(t/T)e^{-it\Omega}dt
=\widehat w\left(\frac{T\Omega}{2\pi}\right)
\ll_{w,A}(1+T|\Omega|)^{-A}.\label{cor:window-decay}\end{equation}
The total absolute coefficient mass in the expanded product is at most a fixed power of $T$. Indeed, the local bound gives $|b_C(p^k)|p^{-k/2}\ll p^{-a_Ck}\le1$, and
\[
\sum_{p^k\le T^{c\sigma_j}}(\log p)|b_C(p^k)|p^{-k/2}
\ll T^{c\sigma_j}\log T.
\]
Combining this with \eqref{cor:spacing} bounds the off-diagonal contribution by
$O(T^\vartheta(\log T)^rT^{-A(1-\vartheta)})$, before the harmless normalization by $\lambda^r$. Taking $A$ sufficiently large proves an error $O(T^{-M})$.

On the diagonal $a=b$, unique factorization of rational integers says that the signed exponents balance at each prime separately. Independent phase integration imposes exactly those same conditions, since
$\mathbb E e^{ih\theta_p}=\mathbf1_{h=0}$. Each diagonal frequency has height average $\int w=1$. This proves \eqref{cor:moment-comparison}.

Apply the same comparison to every subcollection of the variables. The support inequality remains valid. Joint cumulants are finite polynomials in these moments. All the moments involved have polynomial growth by the coefficient estimate above, so errors remain $O(T^{-M})$ after starting with sufficiently large powers of saving. The random vectors belonging to different primes are independent. Additivity of joint cumulants across independent vectors gives \eqref{cor:cumulant-comparison}.
\end{proof}

\subsection{Higher connected correlations}
For a nonnegative integer $a$, define the local function
\begin{equation}Y_{C,p}^{[a]}(\theta)=(\log p)^{a+1}\sum_{k\ge1}\frac{k^a}{p^{k/2}}
\left\{b_C(p^k)e^{ik\theta}+(-1)^a\overline{b_C(p^k)}e^{-ik\theta}\right\}.\label{cor:jet-phase}\end{equation}
This converges at each fixed prime. For $r\ge3$ and a multi-index $\boldsymbol a=(a_1,\ldots,a_r)$ put
\begin{equation}\Gamma_{\boldsymbol C,\boldsymbol a}
=\sum_p\operatorname{cum}_{\theta_p}
  (Y_{C_1,p}^{[a_1]},\ldots,Y_{C_r,p}^{[a_r]}).\label{cor:Gamma}\end{equation}
Theorem~\ref{hi:convergence} proves absolute convergence, including the finitely many bad primes. For $\boldsymbol a=\boldsymbol0$, this is the all-prime version of $\mathcal K(\boldsymbol C;\boldsymbol0)$ in \eqref{hi:global-cumulant}.

\begin{theorem}[Higher connected terms over arbitrary fields]\label{cor:higher}
Let $r\ge3$ be fixed, and assume \eqref{cor:support}. For every integer $J\ge0$,
\begin{align}
\operatorname{cum}_T(X_{C_1}(t;f_1),\ldots,X_{C_r}(t;f_r))
&=(-1)^r\sum_{|\boldsymbol a|\le J}
\frac{\Gamma_{\boldsymbol C,\boldsymbol a}}{\lambda^{r+|\boldsymbol a|}}
\prod_{j=1}^r\frac{\widehat f_j^{(a_j)}(0)}{a_j!}
+O(\lambda^{-r-J-1}).\label{cor:higher-expansion}
\end{align}
This theorem is unconditional for fixed $\Pi/E$ and $\pi/F$. In particular,
\begin{equation}\operatorname{cum}_T(X_{C_1},\ldots,X_{C_r})
=\frac{(-1)^r\Gamma_{\boldsymbol C,\boldsymbol0}}{\lambda^r}
\prod_j\widehat f_j(0)+O(\lambda^{-r-1}).\label{cor:higher-leading}\end{equation}
If all tests are even, the error in \eqref{cor:higher-leading} is $O(\lambda^{-r-2})$.
\end{theorem}
\begin{proof}
Use \eqref{cor:cumulant-comparison}. Expanding its finite partition formula reduces the assertion to balanced prime monomials. Their test weight is
$\prod_j\widehat f_j(\epsilon_jk_j\log p/\lambda)$. Taylor expansion in the single variable $1/\lambda$ through degree $J$ gives
\[
\sum_{|\boldsymbol a|\le J}\lambda^{-|\boldsymbol a|}
\prod_j\frac{\widehat f_j^{(a_j)}(0)}{a_j!}
 (\epsilon_j k_j\log p)^{a_j},
\]
with remainder bounded by
$C\lambda^{-J-1}(\sum_j k_j\log p)^{J+1}$. To justify a bound on the whole range, apply the one-variable Taylor formula with integral remainder to the product; all derivatives of the compactly supported $\widehat f_j$ are bounded. Theorem~\ref{hi:convergence} sums this remainder absolutely over all balanced tuples and primes. The signed factors are precisely those in \eqref{cor:jet-phase}. Multilinearity of cumulants now gives \eqref{cor:higher-expansion}. The comparison error is smaller than every fixed power of $1/\lambda$. For even tests, every first derivative at zero vanishes, so all terms with $|\boldsymbol a|=1$ vanish.
\end{proof}

\subsection{Pair correlations and the required analytic input}
For two real tests define $\phi(u)=\widehat f(u)\widehat g(-u)$ for $u\ge0$. Formula \eqref{cor:cumulant-comparison} with $r=2$ gives
\begin{equation}\operatorname{cum}_T(X_C(t;f),X_D(t;g))
=2\Re\mathcal D_{C,D}(\lambda;\phi)+O(T^{-M}).\label{cor:pair-exact}\end{equation}
Indeed, phase balance forces the two prime powers to be identical; the two possible orientations are complex conjugates. Here the mean terms in the covariance are already covered by the cumulant comparison.

\begin{corollary}[Pair expansion]\label{cor:pair}
Assume \eqref{pa:PNT} and \eqref{cor:support} for the two tests. Then
\begin{align}
\operatorname{cum}_T(X_C(t;f),X_D(t;g))
&=2\Re\left(\rho_{C,D}\int_0^\infty u\widehat f(u)\widehat g(-u)du\right)\notag\\
&\quad+2\Re\sum_{j=0}^{J}
\frac{(-1)^j C_{C,D,j}\phi^{(j)}(0)}{\lambda^{j+2}}
+O(\lambda^{-J-3}).\label{cor:pair-expansion}
\end{align}
In particular, its first arithmetic correction is
\begin{equation}\frac{2\Re C_{C,D,0}}{\lambda^2}\widehat f(0)\widehat g(0).\label{cor:pair-first}
\end{equation}
\end{corollary}
\begin{proof}
Combine \eqref{cor:pair-exact} with Theorem~\ref{pa:smooth}. The values $\widehat f(0)$ and $\widehat g(0)$ are real because the tests are real.
\end{proof}

For a rational-field probe with known base change, the mixed covariance is therefore covered by Section~\ref{id:section}. The self covariance of $\Pi$ is also covered when $E/\mathbb Q$ is Galois, by the conjugate-product identity. In arbitrary fields, knowing the mixed covariance alone does not supply the self covariance. Higher connected terms in Theorem~\ref{cor:higher} do not require this additional input.

\subsection{Raw moments and distinct-index sieving}
For completeness, define the raw smoothed mixed correlation by
\begin{equation}\mathcal R_{\boldsymbol C}(T;\boldsymbol f)
=\mathbb E_T\prod_{j=1}^r Z_{C_j}(t;f_j).\label{cor:raw}
\end{equation}
This is an integrated raw $r$-zero correlation with a separable test function; repeated zero indices are allowed. We give a precise reconstruction procedure rather than identify it without explanation with a more general translation-invariant test class.

Uniformly for $u$ in the compact support of $w$, Stirling's formula in \eqref{cor:center} gives
\begin{equation}M_{C_j}(Tu;f_j)=a_j(u;\lambda)+O(T^{-1/2}),\qquad
a_j(u;\lambda)=\widehat f_j(0)\left\{\frac{\Delta_{C_j}}c+
\frac{\log Q_{C_j}+\Delta_{C_j}\log(u/(2\pi))}{\lambda}\right\}.\label{cor:arch-asymptotic}\end{equation}
The estimate holds with each fixed number of $u$-derivatives. One way to obtain it is to split the gamma integral at distance $T^{1/2}$ from $Tu$. On the near part, its density differs from $\log Q_C+\Delta_C\log(Tu/(2\pi))$ by $O(T^{-1/2})$. On the far part, rapid decay of $f$ gives an arbitrarily large power saving. The finitely many pole terms also decay faster than any fixed power of $T$, by the Paley--Wiener bound and the distance of their fixed real ordinates from $[T,2T]$. The factor $\widehat f(0)/\lambda$ comes from $(2\pi)^{-1}\int g_t=\widehat f(0)/\lambda$.

\begin{lemma}[Absolute diagonal bounds and small height weights]\label{cor:small-weight}
For each fixed label $C$ and $L\ge1$,
\begin{equation}\sum_{p^k\le e^L}\frac{(\log p)^2|b_C(p^k)|^2}{p^k}\ll_C L^2.\label{cor:absolute-pair}\end{equation}
Every absolute diagonal moment of a fixed number $r$ of the unnormalized phase sums in \eqref{cor:weighted-phase} is $O(\lambda^r)$. Consequently, if a smooth height weight $h_T$ is supported in a fixed compact subset of $(1,2)$ and satisfies $\|h_T^{(j)}\|_\infty\ll_j T^{-b}$ for all fixed $j$, then, under the corresponding support condition,
\begin{equation}\frac1T\int h_T(t/T)\prod_{j=1}^rX_{C_j}(t;f_j)dt=O(T^{-b})+O(T^{-M})\label{cor:small-weight-bound}\end{equation}
for every fixed $M$, with the empty product allowed.
\end{lemma}
\begin{proof}
At good primes, weighted Cauchy--Schwarz bounds the summand in \eqref{cor:absolute-pair} by the sum of
$d_C f_v(\log p)^2|a_C(v^h)|^2/(Nv)^h$ over $k=f_vh$. The positive Dirichlet series $(\log L_{E_C}^S(s,C\times\widetilde C))''$ instead has coefficient $h(\log Nv)^2|a_C(v^h)|^2$. Since $f_v\le h f_v^2$, the latter dominates the former up to $d_C$. For $p^k\le e^L$ one has $p^{-k}\le e\,p^{-k(1+1/L)}$. The diagonal Rankin--Selberg function has a simple pole at one, so
\[
\sum_{p^k\le e^L}\frac{(\log p)^2|b_C(p^k)|^2}{p^k}
\le e d_C(\log L_{E_C}^S)''(1+1/L,C\times\widetilde C)+O_C(1)
\ll_C L^2.
\]
The bad primes contribute $O_C(1)$ by their geometric local bounds. Cauchy--Schwarz gives the corresponding absolute mixed pair bound.

On a diagonal frequency tuple, positions using the same prime form a block of size at least two. A block of size two contributes $O(\lambda^2)$ by \eqref{cor:absolute-pair}; a block of size at least three contributes $O(1)$ by Theorem~\ref{hi:convergence}. Dropping the requirement that different blocks use distinct primes only increases an absolute sum. Summing over the finitely many block partitions proves $O(\lambda^{2\lfloor r/2\rfloor})=O(\lambda^r)$. The same diagonal comparison as before, now with $h_T$, multiplies this bound by $O(T^{-b})$ and by $\lambda^{-r}$. The off-diagonal Fourier decay inherits the factor $T^{-b}$ from all the derivative bounds. This proves \eqref{cor:small-weight-bound} without applying Cauchy--Schwarz to a doubled collection of zero statistics.
\end{proof}

\begin{theorem}[A reconstruction formula for raw correlations]\label{cor:raw-theorem}
Assume \eqref{cor:support}, and assume \eqref{pa:PNT} for every pair of labels appearing in the statistic. For each block $B\subseteq\{1,\ldots,r\}$ of size at least two, let $\mathfrak c_B(\lambda)$ denote the expansion in \eqref{cor:pair-expansion} if $|B|=2$, and in \eqref{cor:higher-expansion} if $|B|\ge3$. Then, as an asymptotic series in $\lambda^{-1}$,
\begin{equation}\mathcal R_{\boldsymbol C}(T;\boldsymbol f)
\sim\int_1^2w(u)\sum_{\mathcal P\in\mathfrak P_r}
\left(\prod_{\{j\}\in\mathcal P}a_j(u;\lambda)\right)
\left(\prod_{\substack{B\in\mathcal P\\|B|\ge2}}\mathfrak c_B(\lambda)\right)du.
\label{cor:raw-partitions}\end{equation}
More precisely, for any fixed target order, expand the finite partition expression and retain all powers through that order; the error is bounded by the next inverse power of $\lambda$. Thus the formula computes every finite secondary order within the stated test class.
\end{theorem}
\begin{proof}
Expand $Z_j=M_j+X_j$. In each term replace $M_j(Tu)$ by \eqref{cor:arch-asymptotic}. Keep every replacement remainder inside its smooth height weight, rather than estimate an absolute product of zero statistics. Equation~\eqref{cor:arch-asymptotic} and its differentiated versions make this weight $O(T^{-1/2})$ in every fixed smooth seminorm. Lemma~\ref{cor:small-weight} therefore makes its contribution $O(T^{-1/2})+O(T^{-M})$, which is smaller than every fixed inverse logarithmic power. No doubled Fourier support is used.

Apply the diagonal comparison to each remaining product of $X_j$ with height weight $w(u)\prod a_j(u;\lambda)$. This weight is a finite linear combination of $w(u)(\log u)^j$, $0\le j\le r$, with uniformly bounded coefficients and derivatives. The proof of Theorem~\ref{cor:comparison} therefore applies unchanged; for a weight not of integral one, its diagonal is multiplied by its integral. Conditional on $u$, the resulting phase variables have zero mean and cumulants independent of $u$. The ordinary moment--cumulant identity for the variables $a_j(u;\lambda)-\lambda^{-1}\sum_pY_{j,p}^{(\lambda)}$ yields exactly \eqref{cor:raw-partitions}. The all-order estimates already proved permit finite truncation to any specified order.
\end{proof}

Here is the distinct-index step in this test class. For a fixed color $C$ and a finite set of positions $I$, put $f_B(x)=\prod_{j\in B}f_j(x)$. Inclusion--exclusion on the partition lattice gives
\begin{equation}\sum_{(\rho_j)_{j\in I}}^{*}\prod_{j\in I}
 f_j\left(\frac{\lambda(\gamma_j-t)}{2\pi}\right)
=\sum_{\mathcal P\in\mathfrak P(I)}
\left\{\prod_{B\in\mathcal P}(-1)^{|B|-1}(|B|-1)!\right\}
\prod_{B\in\mathcal P}Z_C(t;f_B).
\label{cor:sieve}\end{equation}
The star means distinct zero \emph{indices}, with multiplicity retained. To verify the coefficient, fix an equality partition of the indices. Summing the displayed M\"obius coefficients over its refinements gives one for the discrete equality partition and zero otherwise. This is the standard combinatorial sieving used in \cite{RS,LiuYe2002}.

For several colors, apply \eqref{cor:sieve} independently within each color and multiply. Zeros belonging to different factors retain their different indices even if their numerical ordinates coincide. Finally, $\widehat f_B$ is the convolution of the $\widehat f_j$ for $j\in B$, so its support is contained in $[-\sum_{j\in B}\sigma_j,\sum_{j\in B}\sigma_j]$. The sum of the support widths is unchanged by merging. Therefore \eqref{cor:support} remains valid, and Theorem~\ref{cor:raw-theorem} applies to every term of the sieve.

\begin{remark}[Scope of the correlation theorem]
The argument proves raw and distinct-index correlations for the stated smoothed separable tests and strict support box, provided the required pair PNTs are available. It also proves all higher connected centered terms for arbitrary fields without those pair PNTs. Passing to larger support regions or to a broader joint Fourier test class requires an additional argument; neither follows merely from the convergence of the arithmetic constants.
\end{remark}

\section{Probe families and further questions}\label{pr:section}
The preceding sections give two kinds of information. The residue of a realized mixed Rankin--Selberg logarithmic derivative detects a global constituent relationship. The regular pair coefficients and higher connected constants give convergent or regularized combinations of Satake traces. The second kind of information is richer, but converting it into a complete global reconstruction requires additional work.

\subsection{Dirichlet characters as probes over an arbitrary field}
Let $\chi$ be a Dirichlet character modulo $q$, viewed through its associated primitive unitary Hecke character of $\mathbb Q$. Omit all rational primes dividing $q$, as well as any fixed bad set when needed. The character
\begin{equation}\chi_E=\chi\circ N_{E/\mathbb Q}\label{pr:normchar}\end{equation}
always exists. Its local value at a good $v\mid p$ is $\chi(p)^{f_v}$. Hence the rational probe identity specializes to
\begin{equation}K_{\Pi,\chi}^{(q)}(s)
=-\frac{(L_E^{(q)})'}{L_E^{(q)}}(s,\Pi\otimes\overline{\chi_E}).\label{pr:charK}\end{equation}
Here $(q)$ means that primes dividing $q$ are omitted; any additional fixed finite correction is understood explicitly. This identity is a twist of a standard automorphic $L$-function over $E$, and requires no nonabelian base change.

For $m>1$, the ranks in the Rankin--Selberg interpretation are different, so there is no pole at $1$ and no shifted overlap pole. The mixed leading coefficient is zero. Nevertheless, $A_{\Pi,\chi}$ is holomorphic near $1$, and its Taylor coefficients are the secondary arithmetic data of Section~\ref{pa:section}. In particular,
\begin{equation}C_{\Pi,\chi,0}
=\left(\frac{(L_E^{(q)})'}{L_E^{(q)}}\right)'(1,\Pi\otimes\overline{\chi_E})
-H_{\Pi,\chi}^{(q)}(1).\label{pr:charC}\end{equation}
This expression uses the same omitted-prime convention in both terms. For $m=1$, any overlap pole is first subtracted as in Proposition~\ref{pa:finite}.

The arithmetic meaning of varying $\chi$ is explicit. If $(a,q)=1$, ordinary character orthogonality gives, initially for $\Re s>1$,
\begin{equation}\frac1{\varphi(q)}\sum_{\chi\bmod q}\chi(a)A_{\Pi,\chi}^{(q)}(s)
=\sum_{\substack{p\nmid q,\ k\ge1\\p^k\equiv a\pmod q}}
\frac{(\log p)^2b_\Pi(p^k)}{p^{ks}}.\label{pr:orthogonality}\end{equation}
Indeed, the character factor on the left is
$\varphi(q)^{-1}\sum_\chi\chi(a)\overline{\chi(p^k)}$, which equals one precisely for $p^k\equiv a\pmod q$. Imprimitive characters cause no problem because every factor at a prime dividing $q$ has been removed. When $m>1$, continuation of the finite left-hand sum defines the regularized right-hand side near $1$. Its Taylor coefficients are the corresponding character combinations of the pair constants.

Thus character probes give congruence-weighted combinations of grouped Satake traces even when every leading mixed overlap vanishes. Equation~\eqref{pr:orthogonality} is not a claim that finitely many such constants determine each individual prime coefficient. That inverse question remains to be studied. Moreover, for real symmetric pair tests one observes real parts of the constants. Recovering complex information requires a suitably varied, oriented family of tests or probes, and must be justified rather than inferred from a single real statistic.

\subsection{Galois symmetry and higher-rank probes}
Suppose $E/\mathbb Q$ is Galois. The integer
\begin{equation}B_\Pi=\#\{g\in\operatorname{Gal}(E/\mathbb Q):\Pi^g\simeq\Pi\}\label{pr:stabilizer}\end{equation}
is visible in the leading self pair term. For a quadratic field it is either one or two. The first case means that $\Pi$ and its nontrivial conjugate are distinct; the second means that $\Pi$ is invariant under the conjugation. These conclusions concern the representation itself and do not require a factorization of its standard $L$-function over $\mathbb Q$.

If a rational-field probe $\pi$ of the same rank has a cuspidal base change, its mixed overlap is one exactly when that base change is $\Pi$. The resulting local equalities are \eqref{in:Satake}, not merely a divisibility condition for Euler polynomials. For a higher-rank probe whose base change has several unitary cuspidal constituents, the residue counts the copies of $\Pi$. The local polynomial of $\Pi$ then divides the base-change polynomial at each good place, but the global pole criterion supplies the compatibility that local divisibility by itself cannot establish; see the symmetric-square example in Section~\ref{id:section}.

\subsection{What zeros can and cannot distinguish}
All the formulas in this paper depend on $\Pi$ through the grouped coefficients $b_\Pi(p^k)$, its archimedean factor, and its conductor. In particular, if $E/\mathbb Q$ is Galois, reindexing places gives
\begin{equation}L_E(s,\Pi)=L_E(s,\Pi^g).\label{pr:conjugate-ambiguity}\end{equation}
Multiplication by an arbitrary fixed probe preserves this equality of standard $L$-functions and their zero multisets. Zero statistics therefore cannot recover place labels that are lost in this reindexing. This is consistent with the automorphic-induction goal: the formally induced polynomial \eqref{in:polynomial} is unchanged by Galois conjugation.

Even complete knowledge of leading overlaps only tests constituent coincidences in the probe families whose Rankin--Selberg realizations are known. The secondary constants distinguish more arithmetic expressions, but no injectivity theorem for the collection of those constants is proved here. A useful next step is to formulate the smallest probe family for which equality of all accessible arithmetic data forces equality of the polynomials $P_{\Pi,p}$ at almost all $p$.

\subsection{A concrete continuation of the project}
There are several specific tasks suggested by the results.
\begin{enumerate}[label=\textup{(\arabic*)},leftmargin=*]
\item Determine analytic continuation and a PNT for the mixed Euler product \eqref{id:mixed-local-factor} beyond the realizations in Section~\ref{id:section}. For $F\subset E$, identify exactly which additional tensor-algebra components remain beyond the direct base change to $E$.
\item Use the character combinations \eqref{pr:orthogonality} and the explicit cubic and quartic constants to formulate testable Satake conditions. Establish a converse criterion for a specified probe family rather than assume that all secondary data are automatically complete.
\item Extend the separable-test calculation of Section~\ref{cor:section} to the desired joint Fourier test class while retaining every secondary term. The strict box proof does not resolve the larger-support problem.
\item Investigate effective numerical truncation for fixed data. The positive convergence margins prove asymptotic tail bounds, but the present fixed-field sieve does not provide practical uniform constants in discriminants or conductors. Numerical certification requires additional quantitative input.
\item Compare the resulting invariants in cases where induction or base change is already known. Such examples can separate genuine additional information from identities already forced by the leading overlap and the conductor.
\end{enumerate}
These tasks leave the main existence question of automorphic induction open. The results above give arithmetic estimates and exact detection mechanisms with which that question can be investigated without assuming its conclusion.

\appendix
\section{A character-extension obstruction and the scope of the repair}\label{app:jiang}

The norm-character construction in Section~\ref{sec:norm-sieve} does not use the ray-class character construction in \cite[\S4.1]{Jiang}. We record a concrete reason for keeping these arguments separate. This concerns the general-number-field construction in the cited version, not the rational-field fixed-data argument.

For a finite prime ideal $\mathfrak p$ of a number field $E$, there is an exact sequence
\begin{equation}\label{ns:ray-sequence}
1\longrightarrow K_{\mathfrak p}\longrightarrow\mathrm{Cl}_{\mathfrak p}(E)
\longrightarrow\mathrm{Cl}(E)\longrightarrow1,
\qquad
K_{\mathfrak p}=(\mathcal O_E/\mathfrak p)^\times/
                \operatorname{im}(\mathcal O_E^\times).
\end{equation}
Dualizing shows that a character of $K_{\mathfrak p}$ extends to the ray class group. It does not show that the extension can be chosen with the same order. Abstractly, the nontrivial character of the subgroup of order two in a cyclic group of order four extends only to characters of order four.

\begin{proposition}\label{ns:ray-counterexample}
Take
\[
E=\Q(\sqrt{-5}),\qquad k=2,\qquad
\mathfrak p=(29,\sqrt{-5}-13).
\]
Then $\mathfrak p$ splits completely in
$E(\zeta_2,\sqrt{\mathcal O_E^\times})=E(i)$ and $2$ divides $|K_{\mathfrak p}|$, but no primitive quadratic ray-class character modulo $\mathfrak p$ exists. Thus the character-existence assertions in \cite[Lemma~4.1 and Proposition~4.3]{Jiang}, in the cited version, do not hold in this generality.
\end{proposition}
\begin{proof}
The ring of integers is $\mathcal O_E=\mathbb Z[\sqrt{-5}]$, and its units are $\{1,-1\}$. Since $13^2\equiv-5\pmod{29}$, the displayed ideal is prime with residue field $\mathbb F_{29}$. Therefore
\begin{equation}\label{ns:ray-kernel-example}
K_{\mathfrak p}\simeq\mathbb F_{29}^{\times}/\{\pm1\},
\qquad |K_{\mathfrak p}|=14.
\end{equation}
As $29\equiv1\pmod4$, the polynomial $X^2+1$ has two distinct roots modulo $\mathfrak p$. The prime $\mathfrak p$ is unramified and splits in $E(i)/E$, establishing the stated splitting condition.

Suppose that $\chi$ were a primitive quadratic ray-class character modulo $\mathfrak p$. Its restriction to $K_{\mathfrak p}$ must be nontrivial, because a character trivial on this kernel factors through the ordinary ideal class group and is unramified. The cyclic group in \eqref{ns:ray-kernel-example} has a unique nontrivial character of order two. Its pullback to $\mathbb F_{29}^{\times}$ is the Legendre character. In particular, on principal ideals prime to $\mathfrak p$,
\begin{equation}\label{ns:principal-character}
\chi((a))=\left(\frac{a}{29}\right)
\qquad(a\in\mathbb Z,\ 29\nmid a).
\end{equation}
The possible inverse convention in the ray-class exact sequence does not change this identity for a quadratic character.

Let $\mathfrak a=(2,1+\sqrt{-5})$. Its norm is two. Its square is generated by
\[
4,\quad 2+2\sqrt{-5},\quad (1+\sqrt{-5})^2=-4+2\sqrt{-5},
\]
so $\mathfrak a^2\subseteq(2)$. Both ideals have norm four, hence $\mathfrak a^2=(2)$. The ideal $\mathfrak a$ is prime to $\mathfrak p$, and therefore
\[
1=\chi(\mathfrak a)^2=\chi(\mathfrak a^2)=\chi((2))
=\left(\frac2{29}\right)=-1,
\]
where the last equality follows from $29\equiv5\pmod8$. This contradiction completes the proof.
\end{proof}

The counterexample concerns the construction of the characters, not the truth of Hypothesis H. Theorem~\ref{ns:norm-power} supplies an independent fixed-data argument using norm pullbacks of Dirichlet characters. Its saving is expressed in terms of $[E:\Q]m^2$ because the conductor of a pulled-back character modulo a split rational prime $q$ has norm $q^{[E:\Q]}$. The argument is sufficient for all the grouped prime-power convergence statements in this paper. It does not establish the stronger discriminant- and conductor-uniform assertions of the cited preprint.

\end{document}